% Hyperbolicity of Generic Surfaces of High Degree in Projective 3-Space
% Jean-Pierre Demailly, Jawher El Goul
% Institut Fourier, Universit\'e de Grenoble I, France
% Plain-TeX file

\magnification=1200
\pretolerance=500 \tolerance=1000 \brokenpenalty=5000
\hoffset=0.5cm
\voffset=1cm
\hsize=12.5cm
\vsize=19cm
\parskip 3pt plus 1pt
\parindent=6.6mm

% Macros

\font\fourteenbf=cmbx10 at 14.4pt
\font\twelvebf=cmbx10 at 12pt
\font\twelvebfit=cmmib10 at 12pt
\font\tenbfit=cmmib10
\font\sevenbfit=cmmib7
\font\eightrm=cmr8
\font\eightbf=cmbx8
\font\eighttt=cmtt8
\font\eightit=cmti8
\font\eightsl=cmsl8
\font\sevenrm=cmr7
\font\sevenbf=cmbx7
\font\sixrm=cmr6
\font\sixbf=cmbx6
\font\fiverm=cmr5
\font\fivebf=cmbx5

\font\tenCal=eusm10
\font\sevenCal=eusm7
\font\fiveCal=eusm5
\newfam\Calfam
  \textfont\Calfam=\tenCal
  \scriptfont\Calfam=\sevenCal
  \scriptscriptfont\Calfam=\fiveCal
\def\Cal{\fam\Calfam\tenCal}

\font\fourteenmsb=msbm10 at 14.4pt

\font\tenmsb=msbm10
\font\eightmsb=msbm8
\font\sevenmsb=msbm7
\font\fivemsb=msbm5
\newfam\msbfam
  \textfont\msbfam=\tenmsb
  \scriptfont\msbfam=\sevenmsb
  \scriptscriptfont\msbfam=\fivemsb
\def\Bbb{\fam\msbfam\tenmsb}

% Chargement des fontes de 8 et 6 points :

\font\eighti=cmmi8
\font\eightsy=cmsy8
\font\sixi=cmmi6
\font\sixsy=cmsy6

\skewchar\eighti='177 \skewchar\sixi='177
\skewchar\eightsy='60 \skewchar\sixsy='60

\catcode`\@=11

\def\pc#1#2|{{\bigf@ntpc #1\penalty\@MM\hskip\z@skip\smallf@ntpc #2}}

\def\tenpoint{%
  \textfont0=\tenrm \scriptfont0=\sevenrm \scriptscriptfont0=\fiverm
  \def\rm{\fam\z@\tenrm}%
  \textfont1=\teni \scriptfont1=\seveni \scriptscriptfont1=\fivei
  \def\oldstyle{\fam\@ne\teni}%
  \textfont2=\tensy \scriptfont2=\sevensy \scriptscriptfont2=\fivesy
  \textfont\itfam=\tenit
  \def\it{\fam\itfam\tenit}%
  \textfont\slfam=\tensl
  \def\sl{\fam\slfam\tensl}%
  \textfont\bffam=\tenbf \scriptfont\bffam=\sevenbf
  \scriptscriptfont\bffam=\fivebf
  \def\bf{\fam\bffam\tenbf}%
  \textfont\ttfam=\tentt
  \def\tt{\fam\ttfam\tentt}%
  \textfont\msbfam=\tenmsb
  \scriptfont\msbfam=\sevenmsb   \scriptscriptfont\msbfam=\fivemsb
  \def\Bbb{\fam\msbfam\tenmsb}%
  \abovedisplayskip=6pt plus 2pt minus 6pt
  \abovedisplayshortskip=0pt plus 3pt
  \belowdisplayskip=6pt plus 2pt minus 6pt
  \belowdisplayshortskip=7pt plus 3pt minus 4pt
  \smallskipamount=3pt plus 1pt minus 1pt
  \medskipamount=6pt plus 2pt minus 2pt
  \bigskipamount=12pt plus 4pt minus 4pt
  \normalbaselineskip=12pt
  \setbox\strutbox=\hbox{\vrule height8.5pt depth3.5pt width0pt}%
  \let\bigf@ntpc=\tenrm \let\smallf@ntpc=\sevenrm
  \normalbaselines\rm}
\def\eightpoint{%
  \let\yearstyle=\eightrm
  \textfont0=\eightrm \scriptfont0=\sixrm \scriptscriptfont0=\fiverm
  \def\rm{\fam\z@\eightrm}%
  \textfont1=\eighti \scriptfont1=\sixi \scriptscriptfont1=\fivei
  \def\oldstyle{\fam\@ne\eighti}%
  \textfont2=\eightsy \scriptfont2=\sixsy \scriptscriptfont2=\fivesy
  \textfont\itfam=\eightit
  \def\it{\fam\itfam\eightit}%
  \textfont\slfam=\eightsl
  \def\sl{\fam\slfam\eightsl}%
  \textfont\bffam=\eightbf \scriptfont\bffam=\sixbf
  \scriptscriptfont\bffam=\fivebf
  \def\bf{\fam\bffam\eightbf}%
  \textfont\ttfam=\eighttt
  \def\tt{\fam\ttfam\eighttt}%
  \textfont\msbfam=\eightmsb
  \scriptfont\msbfam=\sevenmsb   \scriptscriptfont\msbfam=\fivemsb
  \def\Bbb{\fam\msbfam\tenmsb}%
  \abovedisplayskip=9pt plus 2pt minus 6pt
  \abovedisplayshortskip=0pt plus 2pt
  \belowdisplayskip=9pt plus 2pt minus 6pt
  \belowdisplayshortskip=5pt plus 2pt minus 3pt
  \smallskipamount=2pt plus 1pt minus 1pt
  \medskipamount=4pt plus 2pt minus 1pt
  \bigskipamount=9pt plus 3pt minus 3pt
  \normalbaselineskip=9pt
  \setbox\strutbox=\hbox{\vrule height7pt depth2pt width0pt}%
  \let\bigf@ntpc=\eightrm \let\smallf@ntpc=\sixrm
  \normalbaselines\rm}

%permet dans la frappe de retrouver les habitudes de
%la dactylographie francaise qui conseille de mettre 
%un espace avant les signes de ponctuation ; : ! ?
%et d'empecher pourtant une coupure de ligne a cet
%endroit dans le traitement

\newskip\LastSkip
\def\nobreakatskip{\relax\ifhmode\ifdim\lastskip>\z@
  \LastSkip\lastskip\unskip\nobreak\hskip\LastSkip
  \fi\fi}
\catcode`\;=\active
\catcode`\:=\active
\catcode`\!=\active
\catcode`\?=\active
\def;{\nobreakatskip\string;}
\def:{\nobreakatskip\string:}
\def!{\nobreakatskip\string!}
\def?{\nobreakatskip\string?}

\frenchspacing
\tenpoint

\newif\ifpagetitre                \pagetitretrue
\newtoks\hautpagetitre         \hautpagetitre={\hfil}
\newtoks\baspagetitre          \baspagetitre={\hfil}
\newtoks\auteurcourant       \auteurcourant={\hfil}
\newtoks\titrecourant         \titrecourant={\hfil}
\newtoks\baspagedroite       \newtoks\baspagegauche
\newtoks\hautpagegauche       \newtoks\hautpagedroite
\hautpagegauche={\tenbf\folio\hfill\eightrm\the\auteurcourant\hfill }
\hautpagedroite={\hfill\the\titrecourant\hfill\tenbf\folio }  
\baspagedroite={\hfil} \baspagegauche={\hfil}
\headline={\ifpagetitre\the \hautpagetitre\else\ifodd\pageno\the\hautpagedroite
\else\the\hautpagegauche\fi\fi}
\footline={\ifpagetitre\the\baspagetitre\global\pagetitrefalse
\else\ifodd\pageno\the\baspagedroite\else\the\baspagegauche\fi\fi}

\def\appeln@te{}
\def\vfootnote#1{\def\@parameter{#1}\insert\footins\bgroup\eightpoint
  \interlinepenalty\interfootnotelinepenalty
  \splittopskip\ht\strutbox % top baseline for broken footnotes
  \splitmaxdepth\dp\strutbox \floatingpenalty\@MM
  \leftskip\z@skip \rightskip\z@skip
  \ifx\appeln@te\@parameter\indent \else{\noindent #1\ }\fi
  \footstrut\futurelet\next\fo@t}

\def\corners{\def\margehaute{\vbox to \hmargehaute{\hbox to \lpage
{\lefttopcorner\hfil\righttopcorner}\vfil}}
\def\margebasse{\vss\hbox to \lpage{\leftbotcorner\hfil\rightbotcorner}}}

\def\footnoterule{\kern-6\p@
  \hrule width 2truein \kern 5.6\p@} % the \hrule is .4pt high

\catcode`\@=12

\def\pd#1#2 {\pc#1#2| }
\catcode`\@=11
\def\p@int{{\rm .}}
\def\p@intir{\discretionary{\rm .}{}{\rm .\kern.35em---\kern.7em}}
\def\pointir{\afterassignment\pointir@\let\next=}
\def\pointir@{\ifx\next\par\p@int\else\p@intir\fi\next}
\catcode`\@=12

\def\abstract#1{\vbox{\eightpoint \pc ABSTRACT|\pointir #1}}

\newdimen\srdim \srdim=\hsize
\newdimen\irdim \irdim=\hsize
\def\NOSECTREF#1{\noindent\hbox to \srdim{\null\dotfill ???(#1)}}
\def\SECTREF#1{\noindent\hbox to \srdim{\csname REF\romannumeral#1\endcsname}}
\newlinechar=`\^^J
\def\openauxfile{
  \immediate\openin1\jobname.aux
  \ifeof1
  \message{^^JCAUTION\string: you MUST run TeX a second time^^J}
  \let\sectref=\NOSECTREF
  \else
  \input \jobname.aux
  \message{^^JCAUTION\string: if the file has just been modified you may 
    have to run TeX twice^^J}
  \let\sectref=\SECTREF
  \fi
  \message{to get correct page numbers displayed in the table of contents^^J}
  \immediate\openout1=\jobname.aux}

\newcount\numbersection \numbersection=-1
\newcount\numberchapter \numberchapter=0
\def\markpage#1{
      \advance\numbersection by 1
      \immediate\write1{\string\def \string\CH
      \romannumeral\numberchapter
      #1\romannumeral\numbersection \string{%
      \number\pageno \string}}}

\def\titre#1|{{\null\baselineskip=18pt
                           {\fourteenbf
  \textfont\msbfam=\fourteenmsb
  \scriptfont\msbfam=\tenmsb
  \scriptscriptfont\msbfam=\sevenmsb
  \def\Bbb{\fam\msbfam\fourteenmsb}
                           \vskip 2.25ex plus 1ex minus .2ex
                           \leftskip=0pt plus \hsize
                           \rightskip=\leftskip
                           \parfillskip=0pt
                           \noindent #1
                           \par\vskip 2.3ex plus .2ex}}
                           \markpage{T}}

\def\section#1|{{\twelvebf
      \textfont0=\twelvebf  \scriptfont1=\tenbf  \scriptscriptfont1=\sevenbf
      \textfont1=\twelvebfit\scriptfont1=\tenbfit \scriptscriptfont1=\sevenbfit
      \par\penalty -500
      \vskip 3.25ex plus 1ex minus .2ex
      \noindent
      #1\vskip1.5ex plus 0.2ex}
      \markpage{S}}

\def\ssection#1|{{\bf
      \par\penalty -200
      \vskip 3.25ex plus 1ex minus .2ex
      #1\pointir}}

\def\skippage{\vfill\eject \ifodd\pageno \else \vphantom{}\vfill\eject \fi
              \numbersection=-1
              \advance\numberchapter by 1}

\long\def\th#1|#2\finth{\par\vskip5pt\noindent
              {\bf #1}{\sl \pointir #2}\par\vskip 5pt}

\def\ie{{\sl i.e.\ }}

\def\remarque#1|{\par\vskip5pt\noindent{\bf #1}\pointir }
\def\rque#1|{\par\vskip5pt{\sl #1}\pointir}

\def\ieme{\raise 1ex\hbox{\pc{}i\`eme|}}
\def\omini{\raise 1ex\hbox{\pc{}o|}}
\def\emini{\raise 1ex\hbox{\pc{}e|}}
\def\ermini{\raise 1ex\hbox{\pc{}er|}}
\def\\{\mathop{\hbox{\tenmsb r}}\nolimits}

%variables pour la biblio chiffre/lettre
\newif\ifchiffre
\def\chiffre{\chiffretrue}
\chiffre
\newdimen\laenge
\def\lettre#1|{\setbox3=\hbox{#1}\laenge=\wd3\advance\laenge by 3mm
\chiffrefalse}
%reference pour un article :
\def\article#1|#2|#3|#4|#5|#6|#7|%
    {{\ifchiffre\leftskip=7mm\noindent
     \hangindent=2mm\hangafter=1
\llap{[#1]\hskip1.35em}{\bf #2}\pointir {\sl #3}, {\rm #4}, \nobreak{\bf #5} ({\yearstyle #6}), \nobreak #7.\par\else\noindent
\advance\laenge by 4mm \hangindent=\laenge\advance\laenge by -4mm\hangafter=1
\rlap{[#1]}\hskip\laenge{\bf #2}\pointir
{\sl #3}, #4, {\bf #5} ({\yearstyle #6}), #7.\par\fi}}
%reference pour un livre :
\def\livre#1|#2|#3|#4|#5|%
    {{\ifchiffre\leftskip=7mm\noindent
    \hangindent=2mm\hangafter=1
\llap{[#1]\hskip1.35em}{\bf #2}\pointir{\sl #3}, #4, {\yearstyle #5}.\par
\else\noindent
\advance\laenge by 4mm \hangindent=\laenge\advance\laenge by -4mm
\hangafter=1
\rlap{[#1]}\hskip\laenge{\bf #2}\pointir
{\sl #3}, #4, {\yearstyle #5}.\par\fi}}

%reference complementaire :
\def\divers#1|#2|#3|#4|%
    {{\ifchiffre\leftskip=7mm\noindent
    \hangindent=2mm\hangafter=1
     \llap{[#1]\hskip1.35em}{\bf #2}\pointir{\sl #3}, {\rm #4}.\par
\else\noindent
\advance\laenge by 4mm \hangindent=\laenge\advance\laenge by -4mm
\hangafter=1
\rlap{[#1]}\hskip\laenge{\bf #2}\pointir {\sl #3}, {\rm #4}.\par\fi}}
\def\dem{\par\noindent{\sl Proof}\pointir}

\def\carre{\hbox{\font\ppreu=cmsy10\ppreu\char'164\hskip-6.6666pt\char'165}}
\def\finpr{{\ \penalty 500\carre\par\vskip3pt}}

%en-tete papier a lettre Institut Fourier

\def\frac#1/#2{\leavevmode\kern.1em
   \raise.5ex\hbox{\the\scriptfont0 #1}\kern-.1em
      /\kern-.15em\lower.25ex\hbox{\the\scriptfont0 #2}}

\def\buildo#1^#2{\mathrel{\mathop{\null#1}\limits^{#2}}}
\def\buildu#1_#2{\mathrel{\mathop{\null#1}\limits_{#2}}}

%
% La macro \aujour met la date automatiquement i.e. ecrit
% "25 septembre 1986" (Cette macro a ete introduite
% aujour d'hui, qui est le vingt-cinq septembre dix-...)
% la macro \today en est la version anglaise
%
% Correction de l'erreur au niveau de l'implantation de TeX sur la SM90
% concernant l'affectation de la valeur du jour a la variable \day,
%\ifnum\day <0\advance\day by 160\fi
%
\def\aujour{\ifnum\day=1 1\ermini\else\number\day\fi\
\ifcase\month\or janvier\or f\'evrier\or mars\or avril\or mai\or juin\or
juillet\or aout\or septembre\or octobre\or novembre\or d\'ecembre\fi\
\number\year}
\def\today{\ifcase\month\or January \or February \or March \or April\or 
May\or June\or July\or August \or September\or October\or November\or
December\fi\ \number\day , \number\year}
\catcode`\@=11
\newcount\@tempcnta \newcount\@tempcntb 
\def\timeofday{{%
\@tempcnta=\time \divide\@tempcnta by 60 \@tempcntb=\@tempcnta
\multiply\@tempcntb by -60 \advance\@tempcntb by \time
\ifnum\@tempcntb > 9 \number\@tempcnta:\number\@tempcntb
  \else\number\@tempcnta:0\number\@tempcntb\fi}}
\catcode`\@=12

\def\Gr{\mathop{\rm Gr}\nolimits}
\def\Bs{\mathop{\rm Bs}\nolimits}
\def\div{\mathop{\rm div}\nolimits}
\def\mod{\mathop{\rm mod}\nolimits}
\def\coker{\mathop{\rm coker}\nolimits}

\def\lra{\longrightarrow}
\def\ld{,\ldots,}

\def\wt{\widetilde}

\def\bC{{\Bbb C}}
\def\bR{{\Bbb R}}
\def\bP{{\Bbb P}}
\def\bN{{\Bbb N}}
\def\bG{{\Bbb G}}
\def\bZ{{\Bbb Z}}
\def\bQ{{\Bbb Q}}

\def\cF{{\Cal F}}

\def\cO{{\Cal O}}
\def\cS{{\Cal S}}
\def\cT{{\Cal T}}
\def\cX{{\Cal X}}
\def\cY{{\Cal Y}}
\def\cZ{{\Cal Z}}

\def\rk{\mathop{\rm rank}\nolimits}

\def\reg{{\rm reg}}

\def\Pic{\mathop{\rm Pic}\nolimits}

\def\ie{{\rm i.e.}\ }

\def\longhook{\lhook\joinrel\longrightarrow}

%\openauxfile
\titre{Hyperbolicity of Generic Surfaces
of High Degree in Projective 3-Space}|
\auteurcourant{Jean-Pierre Demailly and Jawher El Goul}
\titrecourant{\eightpoint   Hyperbolicity of  Generic Surfaces of
  High Degree  in Projective 3-Space}

\centerline{Jean-Pierre Demailly and Jawher El Goul}
\centerline{{ Institut Fourier, Universit\'e de Grenoble I}}
\centerline{{ 38402 Saint-Martin d'H\`eres, France}}
\vskip 0.2cm

\section 0. Introduction|

The goal of this paper is to study the hyperbolicity of generic
hypersurfaces in projective space. Recall that, by a well-known
criterion due to Brody [Bro78], a compact complex space $X$ is
hyperbolic in the sense of Kobayashi [Ko70] if and only if there is no
nonconstant holomorphic map from $\bC$ to $X$. More than twenty years
ago, Shoshichi Kobayashi proposed the following famous conjecture:
{\sl A~generic $n$-dimensional hypersurface of large enough degree in
$\bP^{n+1}_\bC$ is hyperbolic}. This is of course obvious in the
case of curves: the uniformization theorem shows that a smooth curve
is hyperbolic if and only if it has genus at least $2$, which is the
case if the degree is at least~$4$.

However, the picture is not at all clear in dimension $n\ge 2$.  In
view of results by Zaidenberg [Zai87], the most optimistic lower bound
for the degree of hyperbolic $n$-dimensional hypersurfaces in
$\bP^{n+1}_\bC$ would be $2n+1$ (assuming $n\ge 2$). The hyperbolicity
of $X$ in Kobayashi's analytic setting is expected to be equivalent to
the purely algebraic fact that $X$ does not contain any subvariety 
not of general type (it does imply e.g.\ that $X$ has no rational curve
and no nontrivial image of abelian varieties). L.~Ein has shown 
in [Ein87] that a very generic hypersurface of $\bP^{n+1}_\bC$ of 
degree at least $2n+2$ does not contain any submanifold not of general 
type; a simpler proof has been given later by C.~Voisin [Voi96]. 
The above algebraic property looks however substantially weaker 
than Kobayashi hyperbolicity because it only constrains the geometry
of algebraic subvarieties rather than that of general entire 
transcendental maps.

In the case of a surface~$X$, the optimal degree lower bound for
hyperbolicity is expected to be equal to $5$, which is also precisely
the lowest possible degree for $X$ to be of general type.  In fact,
Green-Griffiths [GG80] have formulated the following much stronger
conjecture: {\sl If $X$ is a variety of general type, every entire
curve $f:\bC\to X$ is algebraically degenerate, and {\rm(optimistic
version of the conjecture)} there is a proper algebraic subset
$Y\subset X$ containing all images of nonconstant entire curves.}
As a (very) generic surface of degree at least $5$ does not contain
rational or elliptic curves by the results of H.~Clemens ([Cl86],
[CKM88]) and G.~Xu [Xu94], it would then follow that such a surface is
hyperbolic. However, almost nothing was known before for the case of
transcendental curves drawn on a (very) generic surface or
hypersurface. Only rather special examples of hyperbolic hypersurfaces
have been constructed in higher dimensions, thanks to a couple of
techniques due to Brody-Green [BG78], Nadel [Na89], Masuda-Noguchi
[MN94], Demailly-El Goul [DEG97] and Siu-Yeung [SY97]. The related
question of complements of curves in $\bP^2$ has perhaps been more
extensively investigated, see Zaidenberg [Zai89, 93],
Dethloff-Schumacher-Wong [DSW92, 94], Siu-Yeung [SY95],
Dethloff-Zaidenberg [DZ95a,b].

Here, we will obtain a confirmation of Kobayashi's conjecture in 
dimension $2$, for the case of surfaces of degree at least $20$. 
Our analysis is based on more general results, which 
also apply to surfaces not necessarily embedded in~$\bP^3$. 
Before presenting them, we introduce some useful terminology.
Let 
$$
f:(\bC,0)\to X
$$ 
be a germ of curve on a surface $X$, expressed as $f=(f_1,f_2)$ in 
suitable local coordinates.  The
notation $E_{k,m}T^\star_X$ stands for the sheaf of ``invariant'' jet
differentials of order $k$ and total degree $m$, which will be defined 
in greater detail in \S$\,$1. For the sake of simplicity, we describe
here the simpler case of jet differentials of order~$2$. A section of 
$E_{2,m}T^\star_X$ is a polynomial differential operator of the form
$$
P(f)=\sum_{\alpha_1+\alpha_2+3j=m}a_{\alpha_1\alpha_2 j}(f)\,
f_1^{\prime \alpha_1}f_2^{\prime \alpha_2}(f_1'f_2''-f_1''f_2')^j
$$
acting on germs of curves. It is clear that $\bigoplus E_{2,m}T^\star_X$ 
is a graded algebra. An {\it algebraic multi-foliation} on a 
surface $X$ is by definition associated with a rank $1$ subsheaf 
$\cF\subset S^mT^\star_X$. Such a subsheaf $\cF$ is generated locally 
by a jet differential of order $1$, i.e.\ a section 
$s\in\Gamma(U,S^mT^\star_X$) of the form
$$
s(z)=\sum_{0\le j\le m}a_j(z_1,z_2)(dz_1)^{m-j}(dz_2)^j,
$$
vanishing at only finitely many points, and such that 
$$s(z)=\prod_{1\le j\le m}(c_{1,j}(z)dz_1+c_{2,j}(z)dz_2)$$ 
factorizes as a product
of generically distinct linear forms. Equivalently, the foliation
is defined by a collection of $m$-lines in $T_{X,z}$ at each
generic point $z$, so that it is associated with a (possibly singular) 
surface $Y\subset P(T_X)$ which is $m$-sheeted over~$X$.
Of course, if $\wt Y$ is a desingularization of $Y$, then
$\wt Y$ carries an associated (possibly singular) foliation, that is, a rank 1
subsheaf of $T^\star_{\wt Y}$. A {\it leaf} of the multi-foliation
on $X$ is just the projection to $X$ of a leaf of the corresponding
foliation on~$\wt Y$. We further introduce the following definition.

\th Definition|
Let $X$ be a nonsingular projective variety of general type. 
We define the {\rm $k$-jet threshold} $\theta_k$ of $X$ to be the infimum
$$
\theta_k=\inf_{m>0}\theta_{k,m} \in\bR,
$$ 
where $\theta_{k,m}$ is the smallest rational number $t/m$ such that there 
is a non zero section in $H^0(X,E_{k,m}T^\star_X\otimes\cO(t\,K_X))$
$($assuming that $t\,K_X$ is an integral divisor, $t\in\bQ)$.
\finth

Since $E_{k,m}T^\star_X\subset E_{k+1,m}T^\star_X$, we have 
of course 
$$
\theta_1\ge\theta_2\ge\ldots\ge\theta_k\ge\ldots.
$$
If $\theta_1<0$, the variety $X$ possesses a lot of $1$-jet differentials,
i.e.\ sections of $H^0(X,S^mT^\star_X)$, and the theory becomes
much easier. The core of the present paper is to investigate the situation
$\theta_1\ge 0$, $\theta_2<0$. It turns out that nonsingular surfaces of 
$\bP^3$ enter in this category when the degree is at least $15$. Degrees
in the range $[5,14]$ would (a priori) only yield $\theta_k<0$ for
values of $k$ at least equal to~$3$, and the situation becomes harder
to study as $k$ increases.

\th Main Theorem|
Let $X$ be a nonsingular surface of general type, and let 
$\theta_k$ be its $k$-jet threshold, $k\ge 1$. Assume that either $\theta_1<0$ 
or that the following three conditions are satisfied$\,:$
\smallskip
\item{\rm(a)} $\theta_1\ge 0$, $\theta_2<0\,;$
\smallskip
\item{\rm(b)} $\Pic(X)=\bZ\,;$
\smallskip
\item{\rm(c)} The Chern numbers of $X$ satisfy~
$\smash{\displaystyle
{c_1^2\over c_2}> {9\over 13+12\,\theta_2}.}
$
\medskip\noindent
Then every nonconstant holomorphic map $f : \bC \to X$ is a
leaf of an algebraic multi-foliation on $X$.
\finth

Our strategy is based on a careful analysis of the geometry of Semple
jet bundles, as proposed in [Dem95]. Following an idea suggested by
Green-Griffiths [GG80], we use Riemann-Roch calculations to prove the
existence of suitable 2-jet differentials of sufficiently large
degree. Actually, it can be shown that the condition $\theta_2<0$ 
always holds true under the assumption $13\,c_1^2-9\,c_2>0$.  Now, any 
$2$-jet differential equation corresponds to a divisor $Z$ in the 
($4$-dimensional) Semple $2$-jet bundle $X_2$. We apply Riemann-Roch 
again on that divisor $Z$ to show that the base locus of $2$-jets is
at most $2$-dimensional -- this is exactly the place where condition
(c) is needed. From this, the existence of the asserted algebraic
multi-foliation follows.

In order to apply the Main Theorem, we still have to check that
conditions (a), (b), (c) are met for a (very) generic surface in $\bP^3$ of
sufficiently high degree.  Here, the terminology ``generic'' (resp.\ 
``very generic'') is used to indicate that the exceptional set is
contained in a finite (resp.\ countable) union of algebraic subsets in
the moduli space of surfaces in $\bP^3\,$.  We prove the following
results (see sections 3 to~6).

\th Proposition 1| Let $X$ be a nonsingular surface of general type such
that $\Pic(X)=\bZ$ and $\theta_1\ge 0$. Then
$$
\theta_2\ge \min\Big(\theta_{2,3}, \theta_{2,4}, \theta_{2,5},
{1\over 2}\theta_1-{1\over 6}\Big).
$$
\finth

\th Proposition 2| Let $X$ be a nonsingular surface of degree $d\ge 5$ in 
$\bP^3$. Then
\smallskip
\item{\rm(a)} $c_1^2=d(d-4)^2$,~~ $c_2=d(d^2-4d+6)$.
\smallskip
\item{\rm(b)} $\Pic(X)=\bZ$ if $X$ is very generic 
$($Noether-Lefschetz theorem$)$.
\smallskip
\item{\rm(c)} $\displaystyle
{1\over d-4}\le \theta_1\le {2\over d-4}$.
\smallskip
\item{\rm(d)} $\theta_2<0$ for $d\ge 15$.
\smallskip
\item{\rm(e)} For a generic surface of degree $d\ge 6$,
$\theta_{2,m}\ge -{1\over 2m}+{2-7/2m\over d-4}$ if $m=3,\,4,\,5$.
\smallskip
\item{\rm(f)} For a very generic surface of degree $d\ge 6$,
$\theta_2\ge -{1\over 6}+{1\over 2(d-4)}$.
\smallskip
\item{\rm(g)} Condition {\rm(c)} of the Main Theorem is met
for a very generic surface of degree $d\ge 21$.
\vskip0pt
\finth

Property (c) of Proposition 2 is verified through an explicit 
calculation of sections, made in \S$\,$5. Property (d) is a consequence 
of the fact that $13\,c_1-9\,c_2>0$ for $d\ge 15$.
In order to check property (e), we rely on an elementary but very 
useful ``proportionality lemma''. We are 
indebted to Mihai Paun for a substantial improvement of the earlier
statement of our proportionality lemma. Let us first observe
that there is a natural filtration on $E_{2,m}T^\star_X$, defined
by the degree $j$ of the monomials $(f_1'f_2''-f_1''f_2')^j$ in $P$,
inducing an exact sequence
$$
0 \lra S^mT^\star_X \lra
E_{2,m}T^\star_X \buildo\lra^\Phi E_{2,m-3}T^\star_X \otimes K_X \to 0
$$

\th Proportionality lemma| Let $X$ be a nonsingular surface of
general type. Then, for all sections
$$
P_i\in H^0(X,E_{2,m_i}T^\star_X\otimes\cO_X(t_i\,K_X))
$$
with $m_i=3,\,4,\,5$ and $t_i\in\bQ_+$, $1+t_1+t_2<(m_1+m_2-3)
\theta_{1,m_1+m_2-3}$, 
the section
$$
\beta_1 P_2-\beta_2P_1 \in
H^0(X,E_{2,m_1+m_2-3}T^\star_X\otimes\cO_X((1+t_1+t_2)\,K_X))
$$
associated with $\beta_i=\Phi(P_i)$ vanishes.
\finth

The proportionality lemma has the very interesting feature that it can
convert a nonvanishing theorem into a generic vanishing theorem$\,$! 
Actually, if one can produce examples of sections $P_1$ for $t_1$ 
sufficiently small, then there cannot exist sections $P_2$ for values
of $t_2$ which are still smaller. The construction of meromorphic 
connections introduced by Nadel [Na89], as it turns out, does produce 
adequate sections $P_1$ with values in $E_{2,3}T^\star\otimes\cO(t_1K_X)$,
$t_1\in{}]-1,0[$, for certain very particular surfaces.

According to recent results of M.~McQuillan (see section 6), the Main
Theorem solves Kobayashi's conjecture in the case of surfaces.

\th Corollary 1| 
A very generic surface $X$ in $\bP^3$ of degree $d\geq 21$  is  
Kobayashi hyperbolic, that is, there is no nonconstant holomorphic map
from $\bC$ to $X$.
\finth

As a consequence of the proof, we also get

\th Corollary 2|
The complement of a very generic curve in $\bP^2$ is hyperbolic and
hyperbolically imbedded for all degrees $d\geq 21$. 
\finth

Our hope is that a suitable generalization of the present techniques to
higher order jets will soon lead to a solution of the Green-Griffiths 
conjecture: every holomorphic map from $\bC$ to a surface of general type
is algebraically degenerate. We would like to thank Gerd Dethloff and
Steven Lu for sharing generously their views on these questions, and
Bernie Shiffman for interesting discussions on related subjects.

\section 1. Semple jet  bundles|

Let $X$ be a complex $n$-dimensional manifold. According to 
\hbox{Green-Griffiths} [GG80], we let $J_k\to X$ be the bundle of $k$-jets
of germs of parametrized curves in~$X$, that is, the set of equivalence
classes of holomorphic maps $f:(\bC,0)\to(X,x)$, with the equivalence
relation $f\sim g$ if and only if all derivatives $f^{(j)}(0)=g^{(j)}(0)$
coincide for $0\le j\le k$, when computed in some local coordinate system
of $X$ near~$x$. The projection map $J_k\to X$ is simply $f\mapsto f(0)$.
Thanks to Taylor's formula, the fiber $J_{k,x}$ can  be identified with the
set of $k$-tuples of vectors $(f'(0)\ld f^{(k)}(0))\in(\bC^n)^k$.
It follows that $J_k$ is a holomorphic fiber bundle with typical fiber
$(\bC^n)^k$ over $X$ (however, $J_k$ is not a vector bundle for $k\ge 2$,
because of the nonlinearity of coordinate changes).
In the terminology of [Dem95], a directed manifold is a pair $(X,V)$ where
$X$ is a complex manifold and $V\subset T_X$ a subbundle. Let $(X,V)$ be a 
complex directed manifold. We define $J_kV\to X$ to be the bundle of $k$-jets
of germs of curves \hbox{$f:(\bC,0)\to X$} which are tangent to $V$, i.e., such that
$f'(t)\in V_{f(t)}$ for all $t$ in a neighborhood of~$0$, together with
the projection map $f\mapsto f(0)$ onto~$X$.
It is easy to check that $J_kV$ is actually a subbundle of~$J_k$.
One of the essential tools used here are the projectivized jet bundles 
$X_k\to X$ introduced in [Dem95]. Let $\bG_k$ be the group of germs of
$k$-jets  biholomorphisms of $(\bC,0)$, that is, the group of germs
of biholomorphic maps
$$
t\mapsto\varphi(t)=a_1t+a_2t^2+\cdots+a_kt^k,\qquad
a_1\in\bC^\star,~a_j\in\bC,~j\ge 2,
$$
in which the composition law is taken modulo terms $t^j$ of degree $j>k$.
The group $\bG_k$ acts on the left on $J_kV$ by reparametrization, 
$(\varphi,f)\mapsto f\circ\varphi$. The bundle $X_k$ can then be seen
as a natural compactification of the quotient of the open subset of 
regular jets $J_kV^\reg\subset J_kV$ by the action of $\bG_k$. We recall
here briefly the basic construction.

To a directed manifold $(X,V)$, one associates inductively a sequence of
directed manifolds $(X_k,V_k)$ as follows. Starting with $(X_0,V_0) = (X,V)$, 
one sets inductively $ X_k= P(V_{k-1})$ [$P(V)$ stands for the projectivized
bundle of lines in the vector bundle $V$], where $ V_k$ is the subbundle of 
$T_{ X_k}$ defined at any point $(x,[v])\in X_k$, $v\in V_{k-1,x}$, by
$$ V_{k,(x,[v])} = \Big\{ \xi \in T_{{ X_k},{(x,[v])} }~;~ {(\pi_k)}_\star\xi
\in \bC\cdot v\Big\},~~ \bC\cdot v\subset V_{k-1,x} \subset T_{X_{k-1},x}~.$$
Here  $\pi_k :  X_k \to X_{k-1}$ denotes the natural projection.
We denote by $\cO_{X_k}(-1)$ the tautological line subbundle of
$\pi^\star_k V_{k-1}$, such that 
$$\cO_{X_k}(-1)_{ (x,[v])} = \bC \cdot v,$$
for all $(x,[v])\in X_k = P(V_{k-1})$. By definition, the bundle $V_k$ fits
in an exact sequence
$$
0\lra T_{X_k/X_{k-1}}\lra V_k\buildo{\lra}^{\pi_{k\star}}
\cO_{X_k}(-1)\lra 0,
$$
and the Euler exact sequence of $T_{X_k/X_{k-1}}$ yields
$$
0\lra\cO_{X_k}\lra \pi_k^\star V_{k-1}\otimes\cO_{X_k}(1)
\lra T_{X_k/X_{k-1}}\lra 0.
$$
From these sequences, we infer
$$\rk V_k = \rk V_{k-1}=\cdots=\rk V=r,\qquad
\dim X_k = n+k(r-1).$$
We say that $(X_k,V_k)$ is the $k$-jet
directed manifold associated with $(X,V)$, and we let 
$$
\pi_{k,j}=\pi_{j+1}\circ\cdots\circ\pi_{k-1}\circ\pi_k:X_k \lra X_j,
$$
be the natural projection.

Now, let \hbox{$f:\Delta_r\to X$} be a nonconstant tangent trajectory
to~$V$. Then $f$  lifts to a well defined and
unique trajectory $f_{[k]}:\Delta_r\to X_k$ of $X_k$ tangent to $V_k$.
Moreover, the derivative $f_{[k-1]}'$ gives rise to a section
$$f_{[k-1]}':T_{\Delta_r}\to f_{[k]}^\star\cO_{X_k}(-1).$$
With any section $\sigma$ of $\cO_{X_k}(m)$, $m\ge 0$, on any open set
$\pi_{k,0}^{-1}(U)$, $U\subset X$, we can associate a holomorphic
differential operator $Q$ of order $k$ acting on $k$-jets of germs of curves
$f:(\bC,0)\to U$ tangent to $V$, by putting
$$
Q(f)(t)=\sigma(f_{[k]}(t))\cdot f_{[k-1]}'(t)^{\otimes m}\in\bC.
$$
In order to understand better this correspondence,
let us use locally a coordinate chart and the associated trivialization 
$T_X\simeq\bC^n$, so that the projection $\bC^n\to\bC^r$
onto the first $r$-coordinates gives rise to admissible coordinates on $V$.
Then $f',\,f'',\ldots,\,f^{(k)}$ are in one to one correspondence with
the $r$-tuples
$$
(f'_1,\ldots,f'_r),\quad (f''_1,\ldots,f''_r),~\ldots~~
(f^{(k)}_1,\ldots,f^{(k)}_r).
$$
\th 1.1. Proposition ([Dem95])|The direct image sheaf
$(\pi_{k,0})_\star\cO_{X_k}(m)$ on $X$ 
coincides with the $($locally free$)$ sheaf $E_{k,m}V^\star $ of  
$k$-jet differentials of weighted degree $m$, that is, by 
definition, the set of germs of polynomial differential operators 
$$
Q(f)=\sum_{\alpha_1\ldots\alpha_k\in\bN^r} 
a_{\alpha_1\ldots\alpha_k}(f)\,(f')^{\alpha_1}(f'')^{\alpha_2}
\cdots (f^{(k)})^{\alpha_k}
$$ 
on $J_kV$ $[$in multi-index notation, 
$(f')^{\alpha_1}=(f'_1)^{\alpha_{1,1}}(f'_2)^{\alpha_{1,2}}\ldots
(f'_r)^{\alpha_{1,r}}\,]$, which are moreover invariant under arbitrary
changes of parametrization: a germ of operator $Q\in E_{k,m}V^\star$ is 
characterized by the condition that, for every germ $f\in J_kV$ and
every germ $\varphi\in\bG_k$,
$$
Q\big(f\circ\varphi)=\varphi^{\prime m}\;Q(f)\circ \varphi.
$$
\finth

Observe that the weighted degree $m$ is taken with respect to weights 
$1$ for $f'$, $2$ for $f''$, etc., thus counts the total numbers 
of ``primes'' in each monomial of the expansion of $Q$.

A basic result, relying on the Ahlfors-Schwarz lemma, is that any 
entire curve $f:\bC\to X$ tangent to $V$ must automatically satisfy all
algebraic differential equations $Q(f)=0$ arising from global 
jet differential operators 
$$
Q\in H^0(X,E_{k,m}V^\star\otimes \cO(-A))
$$ 
which vanish on some ample divisor~$A$. More precisely, we have the
following.

\th 1.2. Theorem ([GG80], [Dem95], [SY97])|Assume that there exist 
integers $k,m>0$ and an ample line bundle $A$ on $X$ such that 
$$H^0(X_k,\cO_{X_k}(m)\otimes
{(\pi_{k,0})}^\star A^{-1})\simeq H^0(X,E_{k,m}V^\star\otimes A^{-1})$$
has nonzero sections $\sigma_1\ld\sigma_N$. Let $Z\subset X_k$ be the
base locus of these sections. Then every entire curve $f:\bC\to X$ tangent
to $V$ is such that $f_{[k]}(\bC)\subset Z$. In other words,
for every global $\bG_k$-invariant polynomial differential operator
$Q$ with values in $A^{-1}$, every entire curve $f$ tangent to $V$  must 
satisfy the algebraic differential equation $Q(f)=0$.
\finth

By definition, a line bundle $L$ is {\it big} if there exists an ample divisor 
$A$ on $X$ such that $L^{\otimes m}\otimes \cO(-A)$ admits a nontrivial
global section when $m$ is large (then there are lots of sections, namely
$h^0(X,L^{\otimes m}\otimes \cO(-A))\gg m^n$ with $n=\dim X$).

As a consequence, Theorem 1.2 can be applied when $\cO_{X_k}(1)$ is big.
In the sequel, we will be concerned only with the ``standard case'' $V=T_X$.

A conjecture by Green-Griffiths and Lang states that every entire curve
drawn on a variety of general type is algebraically degenerate,
i.e.\ contained in a proper algebraic subvariety. In view of this
conjecture and of Theorem~1.2, it is especially
interesting to compute the base locus of the global sections of 
jet differentials, sometimes referred to in the litterature as the 
Green-Griffiths locus of $X$. According to the definition of invariant 
$k$-jets given in [Dem95], we introduce instead the base locus $B_k$ of 
invariant $k$-jets, that is, the intersection 
$$
B_k:=\bigcap_{m>0} B_{k,m} \subset X_k
$$
of the base loci $B_{k,m}$ of all line bundles
$\cO_{X_k}(m)\otimes\pi^\star_{k,0} \cO (-A)$, where $A$ is a given
arbitrary ample divisor over $X$ (clearly, $B_k$ does not depend on
the choice of $A$).  Our hope is that
$$
Y:=\bigcap_{k>0} \pi_{k,0} (B_k)\subset X
$$
can always be shown to be a proper subvariety of $X$.
In the present situation, this will be achieved by lowering down the
dimension of $B_k$ as much as possible. For a surface, we will actually 
show that non vertical components of $B_2$ have dimension at least $2$
under reasonable geometric assumptions on~$X$. In general,
our expectation is that non vertical components of $B_k$ have dimension
at most equal to $\dim X$, whenever $X$ is of general type and $k$ is
large enough.

\section 2. Base locus of 1-jets|

From now on, we suppose that $X$ is a nonsingular surface of general
type (in particular, $X$ must be algebraic, see [BPV84]), and let
$c_1$ and $c_2$ be the Chern classes of $X$.  We first describe some
known facts about surfaces of general type with $c_1^2>c_2$, in
connection with the existence of ``symmetric differentials'', i.e.,
sections in $E_{1,m}T^\star_X=S^mT^\star_X$. Section 3 will be devoted
to refinements of these results in the case of order $2$ jets.

The starting point is Hirzebruch's Riemann-Roch formula [Hi66]
$$\chi (X, S^m T^\star_X ) = {m^3\over6} (c_1^2-c_2)+O(m^2).$$
On the other hand, Serre duality implies
$$h^2(X,S^m T^\star_X)=h^0(X, S^m T_X\otimes K_X).$$
A vanishing theorem due to Bogomolov [Bo79] (see also e.g.\ [Dem95], \S$\,$14)
implies that, on a surface $X$ of general type,
$$h^0 \big( X, S^p T_X\otimes K_X^{\otimes q} \big) =0\quad\hbox{
for all $p,\,q$ such that $p-2q>0$}.$$ 
In particular, $h^0(X, S^m T_X\otimes K_X)=0$ whenever $m\ge 3$
and we get
$$
h^0(P(T_X),\cO_{P(T_X)}(m))=h^0(X, S^m T^\star_X )
\geq \chi(X, S^m T^\star_X )\geq {m^3\over 6}
(c_1^2-c_2)+ O(m^2).
$$  
As a consequence, the line bundle $\cO_{X_1}(1)$ is 
big when $c_1^2 > c_2$, and the base locus
$$B_1 = \bigcap_{m>0} \Bs \big|\cO_{X_1}(m)\otimes
\cO(-A) \big|,
$$ 
(which is equal in this case to the Green-Griffiths locus) 
is a proper algebraic subset of $X_1= P(T_X)$. 

Let $Z$ be an irreducible component of $B_1$ which is a horizontal
surface, i.e.\ such that $\pi_{1,0}(Z)=X$. Then the subbundle
$V_1\subset T_{X_1}$ defines on the desingularization $\wt Z$ of $Z$
an algebraic foliation by curves, such that the tangent bundle to the
leaves is given by $T_Z\cap V_1$ at a general point. Indeed, at any 
regular point $x_1=[v]\in Z$, $v\in T_{X,x}$, at which $\pi_{1,0}$ is 
a local biholomorphism onto $X$, $V_{1,x_1}$ consists of those vectors in
$T_{X_1}$ which project to the line $\bC\,v\subset T_{X,x}$, and
$T_{Z,x_1}\cap V_1$ is the lifting of that line by the isomorphism
$(\pi_{1,0})_\star:T_{Z,x_1}\to T_{X,x}$. 

By Theorem 1.2, for any nonconstant entire curve $f : \bC \to X$, the
curve $f_{[1]}$ must lie in some component $Z$ of $B_1$. If $Z$ is not 
horizontal, i.e.\ if $C=\pi_{1,0}(Z)$ is a curve in $X$, then
$f(\bC)\subset C$. Otherwise, we know by the above that $Z$ carries a
canonical algebraic foliation, and that the image of $f_{[1]}$ lies 
either in the singular set of $Z$ or of the projection $\pi_{1,0}:Z\to X$
(which both consist of at most finitely many curves), or is a leaf of 
the foliation. Combining these observations with a theorem of 
A.~Seidenberg [Se68] on desingularization of analytic foliations on
surfaces, F.~Bogomolov [Bo77] obtained the following finiteness theorem.

\th 2.1. Theorem (Bogomolov)|
There are only finitely many rational and elliptic curves on a surface
of general type with $c_1^2 > c_2$.
\finth

Theorem 2.1 can now be seen (see [M-De78]) as a direct consequence of
the following theorem of J.-P.~Jouanolou [Jo78] on algebraic foliations, 
and of the fact that a surface of general type cannot be ruled or elliptic.

\th 2.2. Theorem (Jouanolou)|
Let $L$ be a subsheaf of the cotangent bundle of a projective manifold
defining an analytic foliation of  codimension $1$. Let
$H$ be the dual distribution of hyperplanes in $T_X$. If there is an
infinite number of hypersurfaces tangent to $H$, then $H$ must be the 
relative tangent sheaf to a meromorphic fibration of $X$ onto a curve.
\finth

The above result of Bogomolov does not give information on transcendental 
curves. As observed by Lu and Yau [LY90], one can say more if the topological 
index ${c_1^2-2c_2}\over 3$ is positive, using the following result of 
Y.~Miyaoka [Mi82] on the almost everywhere ampleness of 
$T^\star_X$. We recall here their proof in order to point out the analogy 
with results of section 3 (see [ScTa85] for the general case of
semi-stable vector bundles). 

First recall that a line bundle $L$ on a projective manifold is called
numerically effective (nef) if the intersection $L\cdot C$ is
nonnegative for all curve $C$ in $X$.  A surface $X$ of general type
is called minimal if its canonical bundle $K_X$ is nef.

\th 2.3. Theorem (Miyaoka)|     
Let $X$ be a minimal surface of general type with $c_1^2-2c_2>0$, then
the restriction $\cO_{X_1}(1)_{|Z}$ is big for every horizontal irreducible 
2-dimensional subvariety  $Z$ of $X_1$.
\finth

\dem The Picard group of $X$ is given by
$$\Pic(X_1)=\Pic(X)\oplus \bZ [u]$$  
where $u:=\cO_{X_1}(1)$, and the cohomology ring $H^\bullet(X_1)$ is
given by
$$
H^\bullet(X_1)=H^\bullet(X)\,[u]/(u^2+(\pi^\star c_1)u+\pi^\star c_2)
$$
[$u$ denoting rather $c_1(\cO_{X_1}(1))$ in that case]. In particular,
$$
u^3=u\cdot\pi^\star(c_1^2-c_2)=c_1^2-c_2,\qquad
u^2\cdot\pi^\star K_X=u\cdot\pi^\star c_1^2=c_1^2.
$$
Let $Z$ be an horizontal irreducible $2$-dimensional subvariety. In 
$\Pic(X_1)$, we have  
$$Z\sim mu-\pi^\star F$$ 
for some $m>0$ and some divisor $F$ on $X$. In order to study
$\cO_{X_1}(1)_{|Z}$, we compute the Hilbert polynomial of
this bundle. The coefficient of the leading term is
$$
{(u_{| Z})}^2=u^2\cdot(mu-\pi^\star F)=m(c_1^2-c_2)+c_1\cdot F,
\leqno(\dagger)
$$
by the above Chern class relations. The main difficulty is to control
the term $c_1\cdot F$. For this, the idea is to use a {\it semi-stability
inequality}. The multiplication morphism by the canonical section of
$\cO(Z)$ defines a sheaf injection $\cO(\pi^\star F)\hookrightarrow
\cO_{X_1}(m)$. By taking the direct images on $X$, we get
       $$\cO(F)\hookrightarrow \pi_\star\cO_{X_1}(m)=S^m T^\star_X .$$
Using the  $K_X$-semi-stability of $T^\star_X$  (see [Yau78] or [Bo79]),
we infer
$$
F\cdot (-c_1)\leq {c_1(S^m T^\star_X )\cdot(-c_1)\over{m+1}}={m\over2} c_1^2.
$$
From $(\dagger)$, we get
$$
{(u_{| Z})}^2\ge {m\over 2} (c_1^2-2\,c_2 )>0,
$$
and Riemann-Roch implies that either $\cO_{X_1}(1)_{|Z}$ or 
$\cO_{X_1}(-1)_{|Z}$ is big. To decide for the sign, we observe that
$K_X$ is big and nef and compute
$$
u_{| Z}\cdot \pi^\star K_X=u\cdot (mu-\pi^\star F)\cdot (-c_1)
=mc_1^2+c_1\cdot F\,;\leqno(\dagger\dagger)$$
from this we get $u_{| Z}\cdot \pi^\star K_X\ge {m\over 2} c_1^2>0$
by the semi-stability inequality. It follows that 
$\cO_{X_1}(1)_{| Z}$ is big. \finpr

By applying the above theorem of Miyaoka to the horizontal components $Z$
of $B_1$, we infer as in Theorem 1.2 that every nonconstant entire curve 
$f : \bC \to X$ is contained in the base locus of  $ \cO_{X_1}(k)\otimes 
\cO(-A)_{| Z}$ for $k$ large, if $A$ is a given ample divisor.
Therefore $f$ is algebraically degenerate.

\remarque 2.4 Remark|
Unfortunately, the ``order 1'' techniques developped in this section
are insufficient to deal with surfaces in $\bP^3$, because in this case 
$$c_1^2=d(d-4)^2<c_2=d(d^2-4d+6).$$
Lemma 3.4 below shows in fact that $H^0(X, S^m T^\star_X)=0$ for all $m>0$.

\section 3. Base locus of 2-jets|

The theory of directed manifolds and Semple jet bundles 
makes it possible to extend the techniques of section 2 to the
case of higher order jets. The existence of suitable algebraic
foliations is provided by the following simple observation, once
sufficient information on the base locus $B_k$ is known.

\th 3.1. Lemma|Let $(X_k,V_k)$ be the bundle of projectivized $k$-jets
associated with a surface $X$ and $V=T_X$. For any irreducible
``horizontal hypersurface'' $Z\subset X_k$  $($\ie such that 
$\pi_{k,k-1} (Z)=X_{k-1})$, the intersection
$T_Z\cap V_k$ defines a distribution of lines on a Zariski open subset
of $Z$, thus inducing a $($possibly singular$)$ $1$-dimensional foliation 
on a desingularization of~$Z$.
\finth

\dem 
We have $\rk V_k =2$ and an exact sequence 
$$0\lra T_{X_k/X_{k-1}}\lra V_k \lra\cO_{X_k}(-1)\lra 0$$ 
which follows directly from the inductive definition of~$V_k$.
Thus the intersection $T_Z\cap V_k$ defines  a
distribution of lines on the Zariski open subset of $Z$ equal to the set
of regular points at which $\pi_{k,k-1}:Z\to X_{k-1}$ is \'etale (at such
points, $V_k$ contains the vertical direction and $T_Z$ does not, thus
$V_k$ and $T_Z$ are transverse).\finpr

For general order $k$, it is hard to get a simple decomposition of
the jet bundles $E_{k,m} T^\star_X$, and thus to calculate their 
Euler characteristic. However, for $k=2$ and $\dim X=2$, it is
observed in [Dem95] that one has the remarkably simple filtration
$$ 
\Gr^\bullet E_{2,m}T^\star_X=\bigoplus_{0\le j\le m/3}
S^{m-3j}T^\star_X\otimes K_X^{\otimes j}.
$$
An elementary interpretation of this filtration consists in writing an
invariant polynomial differential operator as
$$
Q(f)=\sum_{0\le j\le m/3}\quad\sum_{\alpha\in\bN^2,\,|\alpha|=m-3j}
a_{\alpha,j}(f)\,(f')^\alpha(f'\wedge f'')^j
$$
where 
$$
f=(f_1,f_2),\qquad (f')^\alpha=(f'_1)^{\alpha_1}(f'_2)^{\alpha_2},
\qquad f'\wedge f''=f'_1f''_2-f''_1f'_2.
$$

As suggested by Green-Griffiths [GG80], we use the Riemann-Roch
formula to derive an existence criterion for global jet differentials.
A~calculation based on the above filtration of $E_{2,m}T^\star_X$ yields
$$\chi\big(X,E_{2,m}T^\star_X \big)={m^4 \over {648}}(13c_1^2-9c_2) + 
O(m^3 ).$$   
On the other hand,
$$
H^2(X,E_{2,m}\otimes\cO (-A))=
H^0(X,K_X\otimes E_{2,m}T^\star_X\otimes\cO(A))
$$
by Serre duality. Since $K_X\otimes (E_{2,m}T^\star_X)\otimes\cO(A)$
admits a filtration with graded pieces
$$
S^{m-3j}T_X\otimes K_X^{\otimes 1-j}\otimes\cO(A),
$$
and $h^0 \big( X, S^p T_X\otimes K_X^{\otimes q} \big)=0$, $p-2q>0$, 
by Bogomolov's vanishing theorem on the general type surface $X$,
we find
$$h^2 \big( X,E_{2,m}T^\star_X\otimes\cO(-A)) =0$$
for $m$ large. 
In the special case when $X$ is a smooth surface of degree $d$ in 
$\bP_{\bC}^3 $, we take $A= {\cO (1)}_{|X}$. Then we have 
$c_1=(4-d)h$ and $c_2=(d^2-4d+6)h^2$ where $h=c_1(\cO (1)_{|X})$,
$h^2=d$, thus  
$$\chi\big( E_{2,m}T^\star_X \otimes \cO (-A) \big)= d
(4\,d^2-68\,d+154){m^4 \over {648}}+ O(m^3 ).$$
A straightforward computation shows that the leading coefficient
$4\,d^2-68\,d+154$ is positive if $d\geq 15$, and a count of degrees 
implies that the $H^2$ group vanishes whenever $\big((m-3j)+2(j-1)\big)
(d-4)-1>0$ for all $j\le m/3$. For this, it is enough that
$2(m/3-1)(d-4)-1>0$, which is the case for instance
if $d\ge 5$ and $m\geq 5$. Consequently we get the following

\th 3.2. Theorem ([Dem95])|
If $X$ is an algebraic surface of general type and
$A$ an ample line bundle over~$X$, then
$$
h^0(X,E_{2,m}T^\star_X\otimes\cO(-A))\ge
{m^4\over 648}(13\,c_1^2-9\,c_2)-O(m^3).
$$
In particular$\,:$
\smallskip
\item{\rm(a)} If $13\,c_1^2-9\,c_2>0$, then $\theta_2<0$.
\smallskip
\item{\rm(b)} Every smooth surface $X\subset\bP^3$ of degree $d\ge 15$
has $\theta_2<0$.
\finth

We now recall a few basic facts from [Dem95]. As $X_2\to X_1\to X$ is a 
tower of $\bP^1$-bundles over $X$, the Picard group
$\Pic(X_2)=\Pic(X)\oplus\bZ u_1\oplus \bZ u_2$ consists of
all isomorphism classes of line bundles
$$
\pi_{2,1}^\star \cO_{X_1}(a_1)\otimes\cO_{X_2}(a_2)\otimes
\pi_{2,0}^\star L
$$
where $L\in\Pic(X)$. For simplicity of notation, we set
$$
\eqalign{
&u_1=\pi_{2,1}^\star\cO_{X_1}(1),\qquad u_2=\cO_{X_2}(1),\cr
&\cO_{X_2}(a_1,a_2):=\pi_{2,1}^\star \cO_{X_1}(a_1)\otimes\cO_{X_2}(a_2)\cr}
$$
for any pair of integers $(a_1,a_2)\in\bZ^2$. The
canonical injection $\cO_{X_2}(-1) \hookrightarrow \pi^\star_2
V_{1}$ and the exact sequence
$$
0\lra T_{X_{1}/X} \lra V_{1} \buildo \lra^{{({\pi_{1})}_\star}}
\cO_{X_{1}}(-1) \lra 0
$$
yield a canonical line bundle morphism
 $$\cO_{X_2}(-1) \buildo \longhook^{\ (\pi^\star_2)\circ
{(\pi_{1})}_\star } \pi^\star_2\  \cO_{X_{1}}(-1)$$
which admits precisely the hyperplane section $ D_2 := P(T_{X_{1}/X})\subset
X_2=P(V_1)$  as its zero divisor. Hence we find
$\cO_{X_2} (-1) = \pi^\star_2\  \cO_{X_{1}}(-1) \otimes
\cO(-D_2)$ and
$$
\cO_{X_2}(-1,1)\simeq\cO (D_2)
$$ 
is associated with an effective divisor in~$X_2$. 

\th 3.3. Lemma|With respect to the projection $\pi_{2,0}:X_2
\to X$, the weighted line bundle $\cO_{X_2}(a_1,a_2)$ is
\smallskip
\item{\rm(a)} relatively effective if and only if 
$a_1+a_2\ge 0$ and $a_2\ge 0\,;$
\smallskip
\item{\rlap{\rm(a${}'$)}\phantom{\rm(a)}} relatively big if and 
only if $a_1+a_2>0$ and $a_2>0\,;$
\smallskip
\item{\rlap{\rm(b)}\phantom{\rm(a)}} relatively nef if and only if 
$a_1\ge 2 a_2\ge 0\,;$
\smallskip
\item{\rlap{\rm(b${}'$)}\phantom{\rm(a)}} relatively ample if 
and only if $a_1>2 a_2>0.$
\smallskip
\noindent
Moreover, the following properties hold.
\smallskip
\item{\rlap{\rm(c)}\phantom{\rm(a)}} For $m=a_1+a_2\ge 0$, there is an 
injection
$$
{\big( \pi_{2,0} \big) }_\star \big( \cO_{X_2} (a_1,a_2) \big) 
\hookrightarrow E_{2,m} T^\star_X,
$$ 
and the injection is an isomorphism if $a_1-2a_2\le 0$.
\smallskip
\item{\rlap{\rm(d)}\phantom{\rm(a)}} Let $Z\subset X_2$ be an irreducible
divisor such that $Z\ne D_2$. Then in $\Pic(X_2)$ we have
$$
Z\sim a_1u_1+a_2u_2+\pi_{2,0}^\star L,\qquad L\in\Pic(X),
$$
where $a_1\ge 2a_2\ge 0$.
\smallskip
\item{\rlap{\rm(e)}\phantom{\rm(a)}} Let $F\in\Pic(X)$ be any divisor or
line bundle. In $H^\bullet(X_2)=H^\bullet(X)[u_1,u_2]$, we have the 
intersection equalities
$$
\eqalign{
&u_1^4=0,\quad u_1^3u_2=c_1^2-c_2,\quad u_1^2u_2^2=c_2,\quad
u_1u_2^3=c_1^2-3c_2,\quad u_2^4=5c_2-c_1^2,\cr
&u_1^3\cdot F=0,\quad u_1^2u_2\cdot F=-c_1\cdot F,\quad u_1u_2^2\cdot F=0,
\quad u_2^3\cdot F=0.\cr}
$$
\finth

\dem The exact sequence defining $V_1$ shows that $V_1$ has splitting type 
$$
V_{1 | F_1}=\cO(2)\oplus\cO(-1)
$$ 
along the fibers $F_1\simeq\bP^1$ of $X_1\to X$, since $T_{X_1/X | F_1}
=\cO(2)$. Hence the fibers $F_2$ of $X_2\to X$ are Hirzebruch surfaces 
$P(\cO(2)\oplus\cO(-1))\simeq P(\cO\oplus\cO(-3))$ and
$$
\cO_{X_2}(1)_{| F_2}=\cO_{P(\cO(2)\oplus\cO(-1))}(1).
$$ 
It is clear that the condition $a_2\ge 0$ is necessary for 
$\cO_{X_2}(a_1,a_2)_{| F_2}$ to be nef or to have nontrivial sections.
In that case, by taking the direct image by $\pi_{2,1}:X_2\to X_1$,
global sections of $\cO_{X_2}(a_1,a_2)_{| F_2}$ can be
viewed as global sections over $F_1\simeq\bP^1$ of
$$
S^{a_2}\big(\cO(-2)\oplus\cO(1)\big)\otimes\cO(a_1)=\bigoplus_{0\le j\le a_2}
\cO(a_1+a_2-3j).
$$
The extreme terms of the summation are $\cO(a_1+a_2)$ and $\cO(a_1-2a_2)$.
Claims (a)--(b) follow easily from this, and (a)${}'$, (b)${}'$ are also
clear since ``being big'' or ``being ample'' is an open condition in 
$\Pic(X_2)$.
\medskip

\noindent (c) We have $\cO_{X_2}(a_1,a_2)=\cO_{X_2}(m)\otimes\cO(-a_1D_2)$,
thus $\cO_{X_2}(a_1,a_2)\subset \cO_{X_2}(m)$ if $a_1\ge 0$ and
$\cO_{X_2}(a_1,a_2)\supset \cO_{X_2}(m)$ if $a_1\le 0$. In the first
case, it is immediately clear that we get an injection
$$
(\pi_{2,0})_\star\cO(a_1,a_2)\subset (\pi_{2,0})_\star\cO_{X_2}(m)
\buildo\lra^{\simeq} E_{2,m}T^\star_X.
$$
In the second case, we a priori have
$$
(\pi_{2,0})_\star\cO(a_1,a_2)\supset (\pi_{2,0})_\star\cO_{X_2}(m)
\buildo\lra^{\simeq} E_{2,m}T^\star_X,
$$
but the above splitting formula shows that $(\pi_{2,0})_\star\cO(a_1,a_2)$
is already largest possible when $a_1-2a_2\le 0$ (which is the case e.g.\
if $(a_1,a_2)=(0,m)$). Hence we have an isomorphism in that case.
\medskip

\noindent (d) If $a_1<2a_2$, we have an injection
$$
\cO_{X_2}(a_1+1,a_2-1)=\cO_{X_2}(a_1,a_2)\otimes\cO_{X_2}(-D_2)
\subset \cO_{X_2}(a_1,a_2)
$$
which induces the same space of sections over each fibre $F_2$. This shows 
that every divisor $Z$ in the linear system $|\cO_{X_2}(a_1,a_2)\otimes
\pi_{2,0}^\star L|$ contains $D_2$ as an irreducible component, and
therefore cannot be irreducible unless $Z=D_2$.
\medskip

\noindent (e) More general calculations are made in [Dem95]. Our formulas
are easy consequences of the relations $u_1^2+c_1u_1+c_2=0$ and
$u_2^2+c_1(V_1)u_2+c_2(V_1)=0$, where 
$$
c_1(V_1)=c_1+u_1,\qquad c_2(V_1)=c_2-u_1^2=2c_2+c_1u_1.\eqno\carre
$$

Under the condition $13c_1^2-9c_2>0$, Theorem 3.2 shows that the order $2$
base locus $B_2$ is a proper algebraic subset of $X_2$. In order to improve 
Miyaoka's result 2.3, we are going to study the restriction of the line 
bundle $\cO_{X_2}(1)$ to any $3$-dimensional component of $B_2$ (if such
components exist). We get the following

\th 3.4. Proposition|
Let $X$ be a minimal surface of general type. If $c_1^2-{9\over7}c_2>0$,
then the restriction of $\cO_{X_2}(1)$ to every irreducible 3-dimensional 
component $Z$ of $B_2\subset X_2$ which projects onto $X_1$ $($``horizontal
component''$)$ and differs from $D_2$ is big.
\finth

\dem Write
$$
Z\sim a_1 u_1 + a_2 u_2 -\pi_{2,0}^\star F,\qquad (a_1,a_2)\in\bZ^2,\quad
a_1\ge 2a_2>0,
$$
where $F$ is some divisor in $X$. Our strategy is to show that 
$\cO_{X_2}(2,1)_{|Z}$ is big. By Lemma 3.3~(e), we find
$$
(2u_1+u_2)^3\cdot Z = (a_1+a_2)(13\,c_1^2-9\,c_2) + 12\,c_1\cdot F.
\leqno(\dagger{\dagger}\dagger)
$$
Now, the multiplication morphism by the canonical section 
of $\cO(Z)$ defines a sheaf injection
$$
\cO(\pi_{2,0}^\star F)\hookrightarrow\cO_{X_2}(a_1,a_2).
$$
By taking direct images onto $X$, $\cO(F)$ can thus be viewed as a
subsheaf of
$$ 
{\big( \pi_{2,0} \big) }_\star \big( \cO_{X_2} (a_1,a_2) \big) 
\subset E_{2,m} T^\star_X
$$ 
where $m=a_1+a_2$. Looking at the filtration of $E_{2,m}T^\star_X$, we
infer that there is a nontrivial morphism 
$$
\cO(F)\longhook S^{m-3j} T^\star_X\otimes K_X^{\otimes j}
$$
for some $j\leq{m \over 3}$. As in \S$\,$2, the semistability 
inequality implies
$$
F\cdot K_X\leq \Big({m-3j\over2}+j\Big)K_X^2\leq {m \over 2}c_1^2,
\qquad\hbox{thus}\quad -c_1\cdot F\le{m\over 2}c_1^2.
$$
Formula $(\dagger{\dagger}\dagger)$ combined with the assumption 
$7\,c_1^2-9\,c_2>0$ implies
$$(2u_1+u_2)^3\cdot Z\ge m(7\,c_1^2-9\,c_2) > 0.$$
The latter inequality still holds if we replace $\cO_{X_2}(2,1)$ by
$\cO_{X_2}(2+\varepsilon,1)$ with a fixed sufficiently small
positive rational number $\varepsilon$. By  Riemann-Roch, either 
$$
h^0(Z,\cO_{X_2}((2+\varepsilon )p,p)_{| Z})\quad{\rm or}\quad
h^2(Z,\cO_{X_2}((2+\varepsilon )p,p)_{| Z})
$$ 
grows fast as $p$ goes to infinity. We want to exclude the second
possibility. For this, we look at the exact sequence
$$
0\to\cO(-Z)\otimes\cO_{X_2}((2+\varepsilon)p,p)\to
\cO_{X_2}((2+\varepsilon)p,p)\to\cO_Z\otimes \cO_{X_2}((2+\varepsilon)p,p)
\to 0
$$
and take the direct images to $X$ by the Leray spectral sequence of
the fibration $X_2\to X$. As $\cO_{X_2}(2+\varepsilon,1)$ is
relatively ample, the higher $R^q$ sheaves are zero and we see
immediately that
$$
\eqalign{
h^2(Z,\cO_{X_2}((2+\varepsilon)p,p)_{|Z})
&\leq h^2(X_2,\cO_{X_2}((2+\varepsilon)p,p))\cr
&\kern50pt{}+h^3(X_2,\cO_{X_2}(-Z)\otimes\cO_{X_2}((2+\varepsilon)p,p))\cr
&\leq h^2(X,(\pi_{2,0})_\star\cO_{X_2}((2+\varepsilon)p,p)) \cr}
$$
By Bogomolov's vanishing theorem, 
the latter group is zero.
Thus, we obtain that $\cO_{X_2}(2+\varepsilon,1)_{| Z}$ is big,
and this implies that $\cO_{X_2}(1)_{| Z}$ is also big 
because we have a sheaf injection
$$\cO_{X_2}(2+\varepsilon,1)= \cO_{X_2}(3+\varepsilon)\otimes\cO(-(2+\varepsilon)D_2) \longhook \cO_{X_2}(3+\varepsilon)$$
(if necessary, pass to suitable tensor multiples to avoid denominators).
\finpr  

\th 3.5. Corollary|Let $X$ be a surface of general type such that
$c_1^2-{9\over 7}c_2>0$. Then the irreducible components of the
Green-Griffiths locus $B_2\subset X_2$ are of dimension $2$ at most,
except for the divisor $D_2\subset X_2$.
\finth

This corollary is not really convincing, since we already have
sections in \hbox{$H^0(X,S^mT^\star_X\otimes\cO(-A))$} under the weaker
condition $c_1^2-c_2>0$ (a condition which is anyhow
too restrictive to encompass the case of surfaces in $\bP^3$). 
Fortunately, under the additional assumption that the surface has
Picard group $\bZ$, one can get a more precise inequality than the
stability inequality, and that inequality turns out to be sufficient
to treat the case of generic surfaces of sufficiently high degree in
$\bP^3$. 

\section 4. Proof of the Main Theorem|

We assume here throughout that $X$ is a surface of general type such that
$\Pic(X)=\bZ$. Then the canonical bundle $K_X$ is ample, and we have 
$c_1^2>0$, $c_2\ge 0$. Our first goal is to estimate the $2$-jet
threshold of $X$. Consider a non trivial section
$$
\sigma\in H^0(X,E_{2,m}T^\star_X\otimes\cO(t\,K_X)),
\qquad m>0,~~t\in\bQ
$$
and its zero divisor
$$
Z_\sigma=m\,u_2 + t\,\pi_{2,0}^\star K_X\qquad\hbox{in $\Pic(X_2)$}.
$$
Let $Z_\sigma=\sum p_jZ_j$ be the decomposition of $Z_\sigma$ in 
irreducible components. From the equality 
$\Pic(X_2)=\Pic(X)\oplus\bZ u_1\oplus \bZ u_2$ and the 
assumption $\Pic(X)\simeq\bZ$, we find 
$$Z_j\sim a_{1,j} u_1 + a_{2,j} u_2 + t_j\,\pi_{2,0}^\star K_X,$$
for suitable integers $a_{1,j},\ a_{2,j} \in \bZ$ and rational 
numbers $t_j\in \bQ$. By Lemma 3.3, as $Z_j$ is effective,
we must have one of the following three disjoint cases:
\smallskip
\item{$\scriptstyle\bullet$} $(a_{1,j},a_{2,j})=(0,0)$ and 
$Z_j\in\pi_{2,0}^\star\Pic(X)$, $t_j>0\,$;
\smallskip
\item{$\scriptstyle\bullet$} $(a_{1,j},a_{2,j})=(-1,1)$, then $Z_j$ contains
$D_2$, so $Z_j=D_2$ and $t_j=0\,$;
\smallskip
\item{$\scriptstyle\bullet$} $a_{1,j}\ge 2 a_{2,j}\ge 0$ and 
$m_j:=a_{1,j}+a_{2,j}>0$.
\smallskip\noindent
In the third case, we obtain a section
$$\sigma_j\in 
H^0\big(X_2,\cO_{X_2}(m_j)\otimes \pi_{2,0}^\star \cO(t_jK_X)\big)$$
whose divisor is $Z_j+a_{1,j}D_2$. As $m=\sum m_j$ and $t=\sum t_j$,
it is clear that ${t\over m}\ge \min{t_j\over m_j}$, where the minimum is
taken over those sections arising from the third case. It follows that
the $2$-jet threshold can be computed by using only those sections which
correspond to an irreducible divisor (regardless of $D_2$ which is
``negligible'' in this matter). We use the following lemma.

\th 4.1. Lemma|
Let $m=3p+q$, $0\le q\le 2$ a positive integer.
\smallskip
\item{\rm(a)} There are bundle morphisms
$$
E_{2,m}T^\star_X\to E_{2,m-3}T^\star_X\otimes K_X
\to E_{2,m-6}T^\star_X\otimes K_X^2\to\cdots\to
S^qT^\star_X\otimes K_X^p.
$$
\smallskip
\item{\rm(b)} There is a $($nonlinear$\,!)$ discriminant mapping
$$
\Delta:E_{2,m}T^\star_X\to S^{(p-1)(3p+2q)}T^\star_X\otimes K_X^{p(p-1)}.
$$
\smallskip
\item{\rm(c)} If $\theta_1\ge 0$, the $2$-jet 
threshold satisfies
$$
\theta_2\ge \min\Big(\theta_{2,3}, \theta_{2,4}, \theta_{2,5},
{1\over 2}\theta_1 - {1\over 6}\Big).
$$
\finth

\dem (a) is a consequence of the filtration described earlier. In order to
prove (b), we write an element of $E_{2,m}T^\star_X$ in the form
$$
P(f)=\sum_{0\le j\le p}a_j\cdot f^{\prime\,3(p-j)+q}\,W^j
$$
where the $a_j$ is viewed as an element of
$S^{3(p-j)+q}T^\star_X\otimes K_X^j$, and
$$
W=f_1'f_2''-f_1''f_2'\in \Lambda^2T_X=K_X^{-1}.
$$ 
The discriminant $\Delta(P)$ is calculated by interpreting $P$ as a polynomial 
in the indeterminate $W$. The precise formula is
$$
\Delta(P)={1\over a_p}\left|
\matrix{
a_0&a_1&a_2&a_3&\ldots&\!a_{p-1}\!&a_p&0&0&\ldots\cr
0&a_0&a_1&a_2&a_3&\ldots&\!a_{p-1}\!&a_p&0&\ldots\cr
\vdots&&&&&&&&&\vdots\cr
\ldots&0&0&a_0&a_1&a_2&\ldots&\!a_{p-1}\!&a_p&0\cr
\ldots&0&0&0&a_0&a_1&a_2&\ldots&\!a_{p-1}\!&a_p\cr
b_0&b_1&b_2&\ldots&\!b_{p-2}\!&\!b_{p-1}\!&0&0&0&\ldots\cr
0&b_0&b_1&b_2&\ldots&\!b_{p-2}\!&\!b_{p-1}\!&0&0&\ldots\cr
\vdots&&&&&&&&&\vdots\cr
\ldots&0&0&0&b_0&b_1&\ldots&\!b_{p-2}\!&\!b_{p-1}\!&0\cr
\ldots&0&0&0&0&b_0&b_1&\ldots&\!b_{p-2}\!&\!b_{p-1}\!\cr
}
\right|
{
\hbox{$\left.\vbox to 35pt{}\right\} p-1$}
\atop
\hbox{$\left.\vbox to 35pt{}\right\} p$}\hfill
}
$$
where 
$$
{\partial P\over\partial W}=\sum_{0\le j\le p-1}b_jW^j=
\sum_{0\le j\le p-1}(j+1)a_{j+1}W^j
$$
is the derived polynomial. By counting the degrees of all terms 
$a_j$ and $b_j$ as polynomials in $f'$, one sees that
$\Delta(P)$ is a homogeneous polynomial. Its degree is equal to that of
the diagonal term 
$$
{1\over a_p}a_0^{p-1}b_{p-1}^p=\hbox{Const}(a_0a_p)^{p-1},
$$
which lives in $S^{(p-1)(3p+2q)}T^\star_X\otimes K_X^{p(p-1)}$.
Geometrically, if $P$ is a germ of section of $E_{2,m}T^\star_X$,
$(P=0)$ defines a germ of divisor $Z\subset X_2$, and 
$\Delta(P)=0$ is the divisor in $X_1$ along which the projection
$Z\to X_1$ has branched points.

\noindent
(c) By the observations made at the beginning of the section, we can
start with a section in
$$
\sigma\in H^0(X,E_{2,m}T^\star_X\otimes\cO(t\,K_X)),
$$
associated with an irreducible divisor $Z$ in $X_2$ (up to 
some $D_2$ components). If $m=1,2$, we have $E_{2,m}T^\star_X=S^mT^\star_X$,
thus ${t\over m}\ge \theta_1\ge {1\over 2}\theta_1-{1\over 6}$. If
$m=3,4,5$, then
$$
{t\over m}\ge \min(\theta_{2,3},\theta_{2,4},\theta_{2,5})
$$
by definition. If $m=3p+q\ge 6$, $p\ge 2$, we get a non trivial discriminant
section
$$
\Delta\in H^0(X,S^{(p-1)(3p+2q)}T^\star_X\otimes\cO(p(p-1)K_X+(2p-2)t\,K_X)).
$$
Therefore $2(p-1)t+p(p-1)\ge(p-1)(3p+2q)\theta_1$, and this implies
$$
{t\over m}\ge{3p+2q\over 2m}\theta_1-{p\over 2m}\ge 
{1\over 2}\theta_1-{1\over 6}.
$$
Inequality (c) is proved.\finpr

\noindent
{\bf Proof of the Main Theorem.}
If $\theta_1<0$, then $\cO_{X_1}(1)$ is big and we conclude by a direct
application of Theorem~1.2. Assume now that $X$ satisfies assumptions
(a), (b), (c) of the Main Theorem. As $\theta_2<0$ by (a), we have
a non trivial section
$$
\sigma\in H^0(X,E_{2,m}T^\star_X\otimes\cO(t\,K_X)),
\qquad m>0,~~t\in\bQ,~~t<0,
$$
and the discussion made at the beginning of the section shows that
we can assume that $Z_\sigma=Z+a_1D_2$ for some irreducible
divisor $Z$ in $X_2$ such that
$$
Z\sim a_1 u_1+a_2 u_2 +t\pi_{2,0}^\star K_X\qquad\hbox{in $\Pic(X_2)$},
\qquad a_1+a_2=m.
$$
Formula $(\dagger{\dagger}\dagger)$ of section~3 gives
$$(2u_1+u_2)^3\cdot Z = m(13\,c_1^2-9\,c_2) + 12\,t\,c_1^2,$$ 
and by definition of $\theta_2$ we have $t/m\geq\theta_2$, hence
$$
(2u_1+u_2)^3\cdot Z \geq m((13+12\,\theta_2)\,c_1^2 - 9\,c_2)>0
$$
under assumption (c). As in the proof of Proposition 3.4, we conclude that
the restriction $\cO_{X_2}(1)_{| Z}$ is big.  
Consequently, by Theorem~1.2 (or rather, by the proof of Theorem~1.2, 
see [Dem95]), every nonconstant entire curve $f:\bC\to X$ is such 
that $f_{[2]}(\bC)$ is contained in the base locus of 
$\cO_{X_2}(l)\otimes\pi^\star_{2,0} \cO (-A)_{| Z}$ for $l$ large. 
This base locus is at most $2$-dimensional,
and projects onto a proper algebraic subvariety $Y$ of~$X_1$. 
Therefore $f_{[1]}(\bC)$ is contained in $Y$, and the Main Theorem
is proved. \finpr

\section 5. Vanishing of global 2-jet differentials of degree 3, 4, 5|

This section is devoted to the proof of the generic nonexistence of 
certain 2-jet differentials of small degree, as required in condition (d)
of the Main Theorem. We start with the easier and well-known case of
symmetric differentials (see e.g.\ Sakai [Sa78]), which we just 
investigate briefly for the reader's convenience.

\th 5.1. Lemma|Let $X$ be a nonsingular surface of degree $d$
in~$\bP^3$, $m$ a nonnegative integer and $k\in\bZ$. Then
\smallskip
\item{\rm (a)} $H^0(X,S^mT^\star_X\otimes\cO(k))=0$~~ for all
$k\le\min(2m-1,m-2+d)$.
\smallskip
\item{\rm (b)} $H^0(X,S^mT^\star_X\otimes\cO(k))\simeq
H^0(\bP^3,S^mT^\star_{\bP^3}\otimes\cO(k))$~~ for all $k\le m-2+d$.
\smallskip
\item{\rm (c)} For $d\ge 5$, $X$ is of general type and its $1$-jet 
threshold satisfies
$$
{1\over d-4}\le \theta_1\le {2\over d-4},\qquad
\theta_{1,m}\ge{\min(2,1+(d-1)/m)\over d-4}\quad\hbox{for all $m>0$}.
$$
\finth

\dem The Euler exact sequence
$$
0\lra\cO\lra\cO(1)^{\oplus 4}\lra T_{\bP^3}\lra 0
$$
gives an exact sequence
$$
0\to S^mT^\star_{\bP^3}\otimes\cO(k)\to S^m(\cO^{\oplus 4})
\otimes\cO(k-m)\to S^{m-1}(\cO^{\oplus 4})\otimes\cO(k-m+1)\to 0.$$
As $H^q(\bP^3,\cO(p))=0$ for all $q=1,2$ and for $q=0$, $p<0$, we easily
conclude that $H^q(\bP^3, S^mT^\star_{\bP^3}\otimes\cO(k))=0$ in all
cases
$$
q=0,~~k\le 2m-1,\quad\hbox{or}\quad
q=1,~~k\le m-2,\quad\hbox{or}\quad
q=2,~~k\in\bZ.
$$
[$\,$The case $q=0$ is obtained by considering the restriction of sections to
arbitrary lines in $\bP^3$, and by using $T^\star_{\bP^1}=\cO(-2)\,$].
The exact sequence
$$
0\lra\cO_{\bP^3}(-d)\lra\cO_{\bP^3}\lra\cO_X\to 0
$$
twisted by $S^m_{\bP^3}$ then shows that $H^q(X, S^mT^\star_{\bP^3|X}
\otimes\cO(k))=0$ for $q\le 1$ and $k\le m-2$, and that
$$
H^0(X, S^mT^\star_{\bP^3|X}\otimes\cO_X(k)\simeq
H^0(\bP^3,S^mT^\star_{\bP^3}\otimes\cO(k))
$$
for $k\le m-2+d$. Finally, by taking symmetric powers in the dual 
sequence of
$$
0\lra T_X\lra T_{\bP^3|X}\lra\cO_X(d)\lra 0,
$$
we find a sequence
$$
0\lra S^{m-1}T^\star_{\bP^3|X}\otimes \cO_X(-d)\lra S^m 
T^\star_{\bP^3|X}\lra S^mT^\star_X\lra 0,
$$
from which it readily follows that
$H^0(X,S^m T^\star_X\otimes\cO(k))\simeq
H^0(\bP^3,S^m T^\star_{\bP^3}\otimes\cO(k))$ for $k\le m-2+d$. (b) is
proved, and (a) is a special case.

\noindent(c) We have $K_X=\cO_X(d-4)$. Property (a) shows that there
are no nonzero sections in $H^0(X,S^mT^\star_X\otimes\cO(tK_X))$ unless
$t(d-4)\ge \min(2m,m-1+d)$, and this certainly implies $t/m>1/(d-4)$,
whence the lower bound for $\theta_1$. On the other hand, by taking
$m=d-2$ and $k=2m$, we do find a nonzero section in
$H^0(X,S^mT^\star_X\otimes\cO_X(2m))$, whence the upper bound.\finpr

We now turn ourselves to the question of the existence of $2$-jet
differentials of small degree. For this question, it is especially
convenient to use the concept of meromorphic connections, in the
spirit of the work of Y.T.~Siu [Siu87] and A.~Nadel [Na89].  By
definition, a meromorphic connection is an operator acting on
meromorphic vector fields $v=\sum v_i\,\partial/\partial z_i$, $w=\sum
w_i\,\partial/\partial z_i$ which, in any complex coordinate
system $(z_1\ld z_n)$, has the form
$$
\nabla_w v=\sum_{1\le i,k\le n}
\Big(w_i{\partial v_k\over\partial z_i}+\sum_{1\le j\le n}
\Gamma^k_{ij}w_iv_j\Big){\partial\over\partial z_k}
=d_wv+\Gamma\cdot (w,v).
$$
The Christoffel symbols 
$\Gamma =(\Gamma ^k_{ij})_{1\le i,j,k \le n}$ are thus meromorphic functions.
To such a connection, we associate the Wronskian operator
$$
W_{\nabla}(f)=f'\wedge f''_{\nabla},\qquad f''=\nabla_{f'}f',
$$
given explicitly in coordinates by
$$
\eqalign{
W_{\nabla}(f)&=\Big((f_1'f_2''-f_1''f_2')-\Gamma^2_{1,1}f_1^{\prime 3}+
\Gamma^1_{2,2}f_2^{\prime 3}\cr
&\quad{}+(\Gamma^1_{1,1}-\Gamma^2_{1,2}-\Gamma^2_{2,1})f_1^{\prime 2}f_2'
-(\Gamma^2_{2,2}-\Gamma^1_{1,2}-\Gamma^1_{2,1})f_1'f_2^{\prime 2}\Big)
{\partial\over \partial z_1}\wedge {\partial\over \partial z_2}\,.\cr}
$$
If $B$ is the pole divisor of the coefficients $\Gamma^k_{ij}$, 
the Wronskian operator $W_{\nabla}(f)$ takes values
in $\cO(B)\otimes\cO(\Lambda^2 T_X)=\cO(B-K_X)$, thus 
$$
W_\nabla\in H^0(X,E_{2,3}T^\star_X\otimes\cO(B-K_X)).
$$
The relevant type of connections we need are the ``meromorphic
partial projective connections'' introduced in [EG96] and [DEG97].  A
{\it meromorphic partial projective connection} is a section of the
quotient sheaf of the sheaf of meromorphic connections modulo
meromorphic zero order operators of the form $\alpha(w)v+\beta(v)w$.
The Christoffel symbols are thus supposed to be determined only up to
terms of the form
$$\wt\Gamma^k_{ij}-\Gamma^k_{ij}=\alpha_i\delta_{jk}+\beta_j
\delta_{ik}.$$
Adding such terms to $\nabla$ replaces $f''_\nabla$ with 
$f''_\nabla+\alpha(f')f'+\beta(f')f'$, and thus does not change 
the corresponding Wronskian operator. In dimension $2$, a 
meromorphic connection depends on $8$ Christoffel symbols, but a 
partial projective meromorphic connection depends only on $4$ 
Christoffel symbols. Since the Wronskian operator on a surface
also depends only on 4 coefficients, we see in that case that there is
a one-to-one correspondence between partial meromorphic connections 
and Wronskian operators, and more precisely, between partial meromorphic
connections with pole divisor${}\le B$ and Wronskian operators
$W\in H^0(X,E_{2,3}T^\star_X\otimes\cO(B-K_X))$. To make this even more
precise, let us consider the exact sequences
$$
\eqalign{
&0\lra S^3T^\star_X\lra E_{2,3}T^\star_X\buildo\lra^\Phi K_X\lra 0,\cr
0\lra S^3T^\star_X\otimes\cO(&{}B-K_X)\lra E_{2,3}T^\star_X
\otimes\cO(B-K_X)\buildo\lra^\Phi \cO(B)\lra 0.\cr}
$$
To any nonzero section
$$
P\in H^0(X,E_{2,3}T^\star_X\otimes\cO(B-K_X))
$$
corresponds a section $\beta=\Phi(P)\in H^0(X,\cO(B))$ which can be viewed 
as the ``principal symbol'' of $P$ (coefficient of $f_1'f_2''-f_2'f_1''$).
If the symbol $\beta$ is nonzero, we actually get a Wronskian operator
$$
W(f)=\beta(f)^{-1}P(f)
$$
with pole divisor${}\le B$. 

Our next result is a basic proportionality lemma for $2$-jet
differentials of degree $3$, $4$, $5$. We are indebted to Mihai Paun
[Pa99] for the observation that the proportionality lemma also holds
true for degrees $4$ and $5$ (as a consequence, we are now able to get
substantially better degree bounds than in our earlier version of the
manuscript). For a polynomial differential operator $P(f',f'')$ of
total degree $m=3,\,4,\,5$, the exponent of $(f_1'f_2''-f_2'f_1'')^j$
can only take the values $j=0,\,1$, and we thus get an exact sequence
$$
0 \lra S^mT^\star_X \lra
E_{2,m}T^\star_X \buildo\lra^\Phi E_{2,m-3}T^\star_X \otimes K_X \to 0
$$
where $E_{2,m-3}T^\star_X=S^{m-3}T^\star_X$. Explicitly, if
$$
P=\sum_{|\alpha|=m} a_\alpha(f')^\alpha+
\sum_{|\alpha|=m-3}b_\alpha(f')^\alpha(f_1'f_2''-f_2'f_1''),
$$
then $\beta=\Phi(P)=\sum_{|\alpha|=m-3}b_\alpha(f')^\alpha$.

\th 5.2. Proportionality lemma| Let $X$ be a nonsingular surface of
general type. Then, for all sections
$$
P_i\in H^0(X,E_{2,m_i}T^\star_X\otimes\cO_X(t_i\,K_X))
$$
with $m_i=3,\,4,\,5$ and $t_i\in\bQ_+$, $1+t_1+t_2<(m_1+m_2-3)
\theta_{1,m_1+m_2-3}$, 
the section
$$
\beta_1 P_2-\beta_2P_1 \in
H^0(X,E_{2,m_1+m_2-3}T^\star_X\otimes\cO_X((1+t_1+t_2)\,K_X))
$$
associated with $\beta_i=\Phi(P_i)$ vanishes.
\finth

\dem Consider $P=\beta_1P_2-\beta_2P_1$. This is a differential polynomial
operator, and $\Phi(P)=\beta_1\beta_2-\beta_2\beta_1=0$ by construction.
Hence $P$ can be viewed as a section in
$$H^0(X,S^{m_1+m_2-3}T^\star_X\otimes\cO((1+t_1+t_2)K_X)).$$
By definition of $\theta_{1,m}$, this group vanishes if
$$
1+t_1+t_2<(m_1+m_2-3)\theta_{1,m_1+m_2-3}.
\eqno\carre
$$

In particular, a nonzero section of
$H^0(X,E_{2,3}T^\star_X\otimes\cO(tK_X))$ can be viewed as a partial
meromorphic connection with pole divisor $B\le (1+t)K_X$. From this,
we infer

\th 5.3. Corollary|Let $X$ be a nonsingular surface of general type
with \hbox{$\Pic(X)=\bZ$}. Then there exists at most one partial
projective connection $\nabla$ with pole divisor $B<{1\over 2}
(1+3\theta_{1,3})K_X$.\finpr 
\finth

Examples of partial meromorphic connections with low pole orders can
be explicitly constructed by means of Nadel's technique [Na89] (see
also [EG96], [DEG97] and [SY97]). In particular, one can find
examples -- which are however highly nongeneric -- for which the ratio
$B/K_X$ takes more or less random values in the range $]0,1]$.  
By adjusting the choice of $B$ as close as possible to the upper limit 
${1\over 2}(1+3\theta_1)$, we know that the connection must be
unique, and a nonexistence result follows just by taking $B$ slightly
smaller than the upper limit. In this way we obtain

\th 5.4. Proposition|Let $X$ be a generic surface of degree 
$d\ge 6$ in $\bP^3$. Then 
$$
\theta_{2,m}\ge -{1\over 2m}+{2-7/2m\over d-4}\qquad\hbox{for $m=3,\,4,\,5.$}
$$
\finth

\dem Assume that $X$ is a smooth member of a linear system 
of surfaces
$$
X_\lambda=\big\{
\lambda_0s_0(z)+\lambda_1s_1(z)+\lambda_2s_2(z)+\lambda_3s_3(z)=0
\big\}
$$
where $s_0,\,s_1,\,s_2,\,s_3\in\bC[z_0,z_1,z_2,z_3]$ are homogeneous 
polynomials of degree $d$. According Nadel's method [Nad89], we solve 
the linear system
$$
\sum_{0\le k\le 3} \wt\Gamma ^k_{ij} {\partial s_\ell\over \partial
z_k} = {\partial^2 s_\ell\over \partial z_i \partial z_j},\qquad 0\le
i,j,\ell \le 3,
$$
and get in this way a homogeneous meromorphic connection of degree 
$-1$ on $\bC^4$. One can check
that this connection descends to a partial projective meromorphic 
connection $\nabla$ on $\bP^3$ such that $X_\lambda$ is totally geodesic
(see [DEG97]). Let us consider the specific example
$$
X_a=
\big\{z_0^d+z_1^d+z_2^d+z_3^d+a\,z_0^{k_0}z_1^{k_1}z_2^{k_2}z_3^{k_3}=0\big\},
$$
where $k_0,\,k_1,\,k_2,\,k_3\ge 0$ are integers such 
that $\sum k_i=d$. We take in this case
$$
s_0=z_0^{k_0}(z_0^{d-k_0}+a\,z_1^{k_1}z_2^{k_2}z_3^{k_3}),\qquad 
s_i=z_i^d,\quad i=1,2,3.
$$
A short computation shows that $X_a$ is nonsingular if and only if
$a^d\ne (-d)^d\prod k_i^{-k_i}$ and that the pole divisor 
of the connection $\nabla$ is given by
$$
B=\big\{z_0 z_1 z_2 z_3(d\,z_0^{k_1+k_2+k_3}+
ak_0z_1^{k_1}z_2^{k_2}z_3^{k_3})=0\big\}
$$
($B$ is just the zero divisor of the denominator of the rational functions 
expressing solutions $\wt\Gamma_{ij}^k$ of the above linear
system, after these rational functions have been simplified).
In particular, the ratio
$$
{B\over K_{X_a}}={4+k_1+k_2+k_3\over d-4}
$$
can be taken to be ${p\over d-4}$ for any integer $p$ with $4\le p\le d+4$. 
This yields a section $P_1\in H^0(X,E_{2,3}T^\star_X\otimes\cO(t_1K_X))$ with
$t_1={p\over d-4}-1$. We take $p=[{d+3\over 2}]$ so that
$$
{1\over 2}+t_1={3+\varepsilon/2\over d-4}\quad
\hbox{where}\quad \varepsilon=(d+1)\mod 2,~~\varepsilon\in\{0,1\}.
$$ 
The integer $p$ must be at least equal to $4$, thus our choice is 
permitted if $d\ge 6$. We claim that $X=X_a$ has no non 
trivial section in
$$
H^0(X,E_{2,m}T^\star_X\otimes\cO(tK_X)), \qquad m=m_2=3,\,4,\,5
$$
when ${1\over 2}+t<{2m-3-\varepsilon/2\over d-4}$. Indeed,
for $m_1=3$, $m_2=m$ and $t_2=t$, our choices imply
$$
1+t_1+t_2<{2m\over d-4}\le(m_1+m_2-3)\theta_{1,m_1+m_2-3},
$$
as $\theta_{1,m}\ge\theta_{1,5}\ge {2\over d-4}$ for $m\le 5$ and $d\ge 6$.
By Lemma~5.2, any non zero section $P_2\in
H^0(X,E_{2,m}T^\star_X\otimes\cO(tK_X))$ would yield a
meromorphic connection associated with a Wronskian operator
$P_2/\beta_2=P_1/\beta_1$. As $P_1/\beta_1$ is an irreducible
fraction with $\div\beta_1=B$, we conclude that
$\beta_2/\beta_1\in H^0(X,S^{m-3}T^\star_X\otimes
\cO((t_2-t_1)K_X))$ must be holomorphic, hence 
$$
t_2\ge t_1+(m-3)\theta_{1,m-3}\ge t_1+{2m-6\over d-4}.
$$
On the other hand
$$
t_2=t<-{1\over 2}+{2m-3-\varepsilon/2\over d-4}
=t_1+{2m-6-\varepsilon\over d-4},
$$
contradiction. By the Zariski semicontinuity of cohomology, the group
$$
H^0(X,E_{2,m}T^\star_X\otimes\cO(tK_X))
$$
vanishes for a generic surface~$X$, unless
$$
{t\over m}\ge-{1\over 2m}+{2-(3+\varepsilon/2)/m\over d-4}.
$$
Proposition 5.4 is proved.\finpr

\section 6. McQuillan's work on algebraic foliations|

Recently, using Miyaoka's semi-positivity result for cotangent bundles of  
nonuniruled projective varieties [Mi87] and a dynamic diophantine 
approximation, McQuillan [McQ97] derived strong Nevanlinna
Second Main Theorems for holomorphic mappings $f:\bC\to X$
tangent to the leaves of an algebraic foliation.

\th 6.1. Theorem (McQuillan)|
Every parabolic leaf of an algebraic $($multi-$)$ foliation on a surface 
$X$ of general type is algebraically degenerate.
\finth

The assumption $c_1^2>c_2$ guarantees the existence of an algebraic
multi-foliation such that every $f:\bC\to X$ is contained in one of
the leaves. Thus McQuillan's theorem implies

\th 6.2. Corollary (McQuillan)|
If $X$ is a surface of general type with $c_1^2>c_2$, then all entire curves
of $X$ are algebraically degenerate.
\finth

It turns out that McQuillan's proof is rather involved and goes far
beyond the methods presented here (see also M.~Brunella [Bru98] for
an enlightening presentation of McQuillan's main ideas). Since we do
not need the full force of McQuillan's results, we present here
special cases of our $1$-jet and $2$-jet techniques, which should in
principle be quite sufficient to deal with our application (modulo a
formal computational check which will not be handled here).

\th 6.3. Proposition|Let $X$ be a minimal surface of general type,
equipped with an algebraic multi-foliation $\cF\subset S^mT^\star_X$.
Assume that
$$
m(c_1^2-c_2)+c_1\cdot c_1(\cF)>0.
$$
Then there is a curve $\Gamma$ in $X$ such that all parabolic leaves of
$\cF$ are contained in~$\Gamma$.
\finth

\dem
Notice that every rank $1$ torsion free sheaf on a surface is locally free. 
The inclusion morphism of $\cF$ in $S^mT^\star_X$, viewed as a section of
\hbox{$S^mT^\star_X\otimes\cF^{-1}$}, defines a section of
$\cO_{X_1}(m)\otimes\pi^\star\cF^{-1}$ whose zero divisor $Z\subset 
X_1=P(T_X)$ is precisely the divisor associated with the foliation 
(as explained in the introduction). Therefore $Z=mu-\pi^\star\cF$ in
$\Pic(X_1)$, and our calculations of section~$2$ (see~$(\dagger)$ and
$(\dagger\dagger)$) imply that $\cO_{X_1}(1)_{|Z}$ is big as soon as
$$
{(u_{| Z})}^2=m(c_1^2-c_2)+c_1\cdot \cF>0,\qquad
{(u_{| Z})}\cdot(-c_1)=mc_1^2+c_1\cdot \cF>0.
$$
However, as $X$ is minimal, we have $c_2\ge 0$, and Proposition 6.3
follows.\finpr

Again, the above $1$-jet result is not sufficient to cover the case
of surfaces in~$\bP^3$, so we have to deal with a $2$-jet version
instead. Let $Z\subset X_1=P(T_X)$ be the divisor associated with the 
given foliation~$\cF$, and $\sigma\in H^0(X_1,\cO_{X_1}(m)\otimes
\pi^\star\cF^{-1})$ the corresponding section. We let $\cT_Z$ be
the tangent sheaf to $Z$, i.e.\ the rank $2$ sheaf $\cT_Z$ defined by the 
exact sequence
$$
0\lra\cT_Z\lra T_{X_1|Z} 
\buildo\lra^{d\sigma} \cO_{X_1}(m)_{|Z}\otimes\pi^\star\cF^{-1}_{|Z}
\lra 0.
$$
If we define $\cS=\cT_Z\cap\cO(V_1)$ sheaf-theoretically, we find an
exact sequence
$$
0\lra\cS\lra V_{1|Z} \buildo\lra^{d\sigma}
\cO_{X_1}(m)_{|Z}\otimes\pi^\star\cF^{-1}_{|Z},
$$
where $\cS$ is an invertible subsheaf, and a dual exact sequence
$$
0\lra \cO_{X_1}(-m)_{|Z}\otimes\pi^\star\cF_{|Z}\lra 
V^\star_{1|Z}\lra \cS^\star.
$$
We can then lift $Z$ into a surface $\wt Z\subset X_2$, in such a
way that the projection map $\pi_{2,1}:\wt Z\to Z$ is a modification;
at a generic point $x\in Z$, the point of $\wt Z$ lying above
$x$ is taken to be $(x,[\cS_x])\in X_2$. Our goal is to compute the
cohomology class of the $2$-cycle $\wt Z$ in $H^\bullet(X_2)$.
One of the difficulties is that the cokernel of the map 
$$
d\sigma_{|V_{1|Z}}:V_{1|Z}\to\cO_{X_1}(m)_{|Z}\otimes\pi^\star\cF^{-1}_{|Z}
$$
may have torsion along a $1$-cycle $G_1\subset Z$, i.e., there is a
factorization
$$
d\sigma_{|V_{1|Z}}:V_{1|Z}\to\cO_{X_1}(m)_{|Z}\otimes\pi^\star\cF^{-1}_{|Z}
\otimes\cO_Z(-G_1)\to\cO_{X_1}(m)_{|Z}\otimes\pi^\star\cF^{-1}_{|Z}
$$
such that the cokernel of the first arrow has $0$-dimensional support 
(of course, $G_1$ need not be reduced). If the foliation is generic, 
however, the cokernel of $d\sigma_{|V_{1|Z}}$ will have no torsion
in codimension~$1$, and $d\sigma$ then induces a section of 
$$
\cO_{X_2}(1)\otimes\pi_{2,1}^\star\cO_{X_1}(m)\otimes
\pi_{2,0}^\star\cF^{-1}_{|\pi_{2,1}^{-1}(Z)}\sim
(u_2+m\,u_1-\cF)_{|\pi_{2,1}^{-1}(Z)}
$$
whose zero locus is~$\wt Z$. As $Z\sim mu_1-\cF$, the cohomology class of 
$\wt Z$ in $H^4(X_2)$ is given by
$$
\eqalign{
\{\wt Z\}&=(mu_1-\cF)\cdot(u_2+mu_1-\cF)\cr
&=m^2u_1^2+m\,u_1\cdot u_2-2m\,u_1\cdot\cF-u_2\cdot\cF+\cF^2.\cr}
$$
A short Chern class computation yields
$$
(2u_1+u_2)^2\cdot\wt Z=m^2(4c_1^2-3c_2)+m(5c_1^2-3c_2)+(8m+4)c_1\cdot\cF+3\,
\cF^2.
$$
If the $1$-cycle $G_1$ is nonzero, our numerical formula for $\wt Z$ 
becomes 
$$
\{\wt Z\}=(mu_1-\cF)\cdot(u_2+mu_1-\cF)-\pi_{2,1}^\star\{G_1\}.
$$
On the other hand, we find
$$
(2u_1+u_2)^2\cdot \pi_{2,1}^\star\{G_1\}=(3u_1-c_1)\cdot G_1.
$$
The general formula for $(2u_1+u_2)^2\cdot\wt Z$ is thus
$$
(2u_1+u_2)^2\cdot\wt Z=m^2(4c_1^2-3c_2)+m(5c_1^2-3c_2)+(8m+4)c_1\cdot\cF+3\,
\cF^2-(3u_1-c_1)\cdot G_1.
$$
By using obvious exact sequences, $H^2(\wt Z,m\cO_{X_2}(2,1)_{|\wt Z})$ is
a quotient of
$$
H^2(\pi_{2,1}^{-1}(Z),m\cO_{X_2}(2,1)_{|\pi_{2,1}^{-1}(Z)}),
$$
which is itself controlled by $H^2(X_2,m\cO_{X_2}(2,1))$, \hbox{$H^3
(X_2,m\cO_{X_2}(2,1)\otimes\cO(-Z))$}. A direct image argument shows that
the latter groups are controlled by groups of the form
$H^2(X,E_{2,3m}T^\star_X\otimes L)$, with suitable line bundles $L$.
As in the proof of Theorem~3.4 (possibly after changing $\cO_{X_2}(2,1)$ into
$\cO_{X_2}(2+\varepsilon,1)$ in the above arguments), one can check that the
latter $H^2$ groups vanish. The positivity of $(2u_1+u_2)^2\cdot\wt Z$ thus
implies that $\cO_{X_2}(2,1)_{|\wt Z}$ is big, and therefore all 
parabolic leaves of the (multi)-foliation $\cF$ are algebraically degenerate.
We thus obtain:

\th 6.4. Proposition|Let $X$ be a surface of general type, equipped with
a multi-foliation $\cF\subset S^mT^\star_X$, and let 
$\sigma\in H^0(X_1,\cO_{X_1}(m)\otimes\pi_{1,0}^\star\cF^\star)$ be the 
associated canonical section. Finally, let $G_1$ be the divisorial part of the
subscheme defined by $\coker(d\sigma_{|V_{1|Z}})$. Then, under the assumption
$$
m^2(4c_1^2-3c_2)+m(5c_1^2-3c_2)+(8m+4)c_1\cdot\cF+3\,\cF^2
-(3u_1-c_1)\cdot G_1>0,
$$
all parabolic leaves of $\cF$ are algebraically degenerate.
\finth

\th 6.5. Corollary|Let $X\subset\bP^3$ be a surface of degree $d\ge 18$
with $\Pic(X)=\bZ$, and let $\cF\subset S^mT^\star_X$ be a generic
multi-foliation, in the sense that the $1$-cycle $G_1$ defined above
is zero. Then all parabolic leaves of $\cF$ are algebraically degenerate
and contained in a fixed $1$-dimensional algebraic subset $Y\subset X$.
\finth

\dem Note that the line
subbundle $\cF\subset S^mT^\star_X$ must be negative (otherwise $\cF$
would yield a nontrivial section of $S^mT^\star_X$), hence
$c_1\cdot\cF>0$, $\cF^2>0$, and likewise we have
$$
4\,c_1^2-3\,c_2=d(d^2-20\,d+46)>0,\qquad
5\,c_1^2-3\,c_2=d(2\,d^2-28\,d+62)>0
$$
for $d\ge 18$. Thus, we get the conclusion if $G_1=0$ (but a rather large 
additional contribution of $G_1$ would still be allowable; we do not
know how much of it can actually occur).\finpr

\section 7. Proof of the Corollaries|

\noindent
{\bf Proof of Corollary 1.} 
\smallskip
Recall that by the Noether-Lefschetz theorem, a very generic surface
$X$ in $\bP^3$ is such that $\Pic(X)= \bZ$, with generator
$\cO_X(1)$.  On the other hand, improving a result of
H.~Clemens ([Cl86] and [CKM88]), G.~Xu [Xu94] has shown that the genus
of every curve contained in a very generic surface of degree $d\ge 5$ 
satisfies the bound $g\ge d(d-3)/2-2$ (this bound is sharp). In particular, 
such a surface does not contain rational or elliptic curves. Now
take a very generic surface in $\bP^3$ of degree $d\geq 21$,
which has no rational or elliptic curves, and such that the conclusions
of the Main Theorem apply, i.e.\ every nonconstant entire curve 
$f : \bC \to X$ is such that $f_{[1]}(\bC)$ lies in the leaf of an
algebraic foliation on a surface $Z\subset X_1$.
Then, by McQuillan's result, $f$ must be algebraically
degenerate. The closure $\Gamma=\overline{f(\bC)}$ would then be an
algebraic curve of genus $0$ or $1$, contradiction.

\remarque 7.1. Remark|If one would like to avoid any appeal to McQuillan's
deep result, it would remain to check on an example that the
multi-foliation defined by $Z$ satisfies the sufficient condition
described in Proposition~6.4. This might require for instance a computer
check, and is likely to hold without much restriction. 

\remarque 7.2. Remark|It is extremely likely that Corollary 1 holds true
for generic surfaces and not only for very generic ones. In fact,
since we have a smooth family of nonsingular surfaces $\cX\to M_d
\subset\bP^{N_d}$ in
each degree $d$, the Riemann-Roch calculations of sections 3, 4 hold
true in the relative situation, and thus produce an algebraic family
of divisors $\cZ_t \subset(\cX_t)_2$ on some Zariski open subset
$M'_d\subset M$, $t\in M'_d$.  By shrinking $M'_d$, we can assume that
all $\cZ_t$ are irreducible, and that we have a flat family $\cZ\to
M'_d$. By relative Riemann-Roch again, we get a family of divisors
$\cY_t\subset \cZ_t$, and thus a family of foliations $\cF_t$ on the
$1$-jet bundles $(\cX_t)_1$. Finally, if Proposition 6.4 can be applied
to these foliations (and we strongly expect that this is indeed the
case), we get an algebraic family of curves $\Gamma_t\subset\cX_t$ such
that all holomorphic maps $f:\bC\to \cX_t$ are contained in~$\Gamma_t$. 
As the degree is bounded, a trivial
Hilbert scheme argument implies that the set of $t$'s for which one of
the components of $\Gamma_t$ is rational or elliptic is closed
algebraic and nowhere dense. Our claim follows.  \finpr
\bigskip

\noindent
{\bf Proof of Corollary 2.} 
\smallskip

Let $C=\sigma^{-1}(0)$ be a nonsingular curve of degree $d$ in  $\bP^2$. 
Consider the cyclic covering $X_C=\{z_3^d=\sigma(z_0,z_1,z_2)\}\to\bP^2$ 
of degree $d$, ramified along~$C$. Then $X_C$ is a nonsingular surface 
in~$\bP^3$, and as $\bC$ is simply connected,
every holomorphic map $f : \bC \to  \bP^2 \\ C$ can be lifted
to $X_C$. It is known that $\Pic(X_C)=\bZ$ for generic~$C\,$; see e.g.\ 
J.~Esser's PhD Thesis [Ess93] (we express our thanks to K.~Amerik and
E.~Viehweg for pointing out the reference to us; see Hartshorne 
[Ha75] for the following related well known fact: if $(X_t)_{t\in\bP^1}$
is a Lefschetz pencil of surfaces on a $3$-fold~$W$ and
$H^0(X_t,K_{X_t})\ne 0$, then $\Pic(X_t)\simeq \Pic(W)$ for
generic~$t$). The nonexistence theorem proved in section 6 also holds
true for at least one~$X_C$, for example
$$
\eqalign{
C&=\big\{z_0^d+z_1^d+z_2^d+a\,z_0^{k_0}z_1^{k_1}z_2^{k_2}=0\big\},\cr
X_C&=\big\{z_0^d+z_1^d+z_2^d+z_3^d+a\,z_0^{k_0}z_1^{k_1}z_2^{k_2}=0\big\}.
\cr}
$$
We then conclude as above that $X_C$ is hyperbolic for generic $C$.
This implies in particular that $\bP^2 \\ C$ is hyperbolic and
hyperbolically embedded in $\bP^2$ (see Green [Gr77]).\finpr

\titre References|

{\eightpoint
\parskip=4pt plus 1pt minus 1pt
\lettre M-De78|

\livre BPV84|Barth W., Peters C., Van de Ven A|Compact complex
surfaces|Ergebrnisse der Mathematik und ihrer Grenzgebiete, 3, Folge,
Band 4, Springer|1984|

\article Bo77|Bogomolov F.A|Families of curves on a surface of general
type|Soviet.\ Math.\ Dokl.|18|1977|1294--1297|

\article Bo79|Bogomolov F.A|Holomorphic tensors and vector bundles on 
projective varieties|Math.\ USSR Izvestija|13|1979|499-555|
 
\article Bro78|Brody, R|Compact manifolds and hyperbolicity|Trans.\
Amer.\ Math.\ Soc.|235|1978|213--219|

\article Bru98|Brunella M|Courbes enti\`eres et feuilletages
holomorphes|L'Enseignement Ma\-th\'e\-matique|45|1999|195--216|

\livre CKM88|Clemens H., Koll\'ar J., Mori S|Higher dimensional
complex geometry|Ast\'erisque {\bf 166}|1988|

\article Cl86|Clemens H|Curves on generic
hypersurface|Ann.\ Sci.\ \'Ec.\ Norm.\ Sup.|19|1986|629--636|

\livre Dem95|Demailly J.-P|Algebraic criteria for Kobayashi
hyperbolic projective varieties and jet differentials|Proceedings of
the AMS Summer Institute on Alg.\ Geom.\ held at Santa Cruz,
ed.\ J.~Koll\'ar, July {\yearstyle 1995}|{\rm to appear}|

\livre DEG96|Demailly J.-P., El Goul J|Connexions m\'eromorphes
projectives et vari\'et\'es alg\'ebriques
hyperboliques|C.~R.\ Acad.\ Sci.\ Paris, \`a para\^{\i}tre|{\rm (}1997{\rm )}|

\article DSW92|Dethloff, G., Schumacher, G., Wong, P.M|Hyperbolicity
of the complement of plane algebraic curves|Amer.\ J.\ Math|117|1995|573--599|

\article DSW94|Dethloff, G., Schumacher, G., Wong, P.M|On the
hyperbolicity of the complements of curves in Algebraic surfaces:
the three component case|Duke Math.\ Math.|78|1995|193--212|

\divers DZ95a|Dethloff, G., Zaidenberg, M|Plane curves with C-hyperbolic
complements|Pr\'epublication de l'Institut Fourier, 299, Grenoble
1995, 44$\,$p, Duke Math.\ E-print alg-geom/9501007
(\`a para\^{\i}tre dans Annali Sc.\ Norm.\ Super.\ Pisa)|

\divers DZ95b|Dethloff, G., Zaidenberg, M|Examples of plane curves of
low degrees with hyperbolic and C-hyperbolic complements|
Proc.\ Conf.\ Geometric Complex Analysis, Hayama, 20--29 March, 1995,
J.~Noguchi, H.~Fujimoto, J.~Kajiwara, T.~Ohsawa eds., World
Scientific, Singapore, 17$\,$p (\`a para\^{\i}tre)|

\article Ein88|Ein, L|Subvarieties of generic complete intersections|Invent.\ 
Math.|94|1988|163--169|

\article EG96|El Goul J|Algebraic families of smooth hyperbolic
surfaces of low degree in $\bP^3_\bC$|Manuscripta
Math.|90|1996|521--532|

\divers Ess93|Esser J|Noether-Lefschetz-Theoreme fur zyklische 
belagerungen|PhD Thesis in Mathematik und Informatik der Universit\"at
Gesamthochschule Essen, February 1993|

\article Gr77|Green M|The hyperbolicity of the complement of $2n+1$
hyperplanes in general position in ${\bf P}^n_{\bf C}$ and related
results|Proc. Amer. Math. Soc.|66|1977|109--113|

\article GG80|Green M., Griffiths P|Two applications of algebraic
geometry to entire holomorphic mappings|The Chen Symposium 1979,
Proc.\ Inter.\ Sympos.\ Berkeley, CA, 1979, Springer-Verlag|{\rm, New
York}|1980|41--74|

\article Ha75|Hartshorne R| Equivalence relations on algebraic cycles
and subvarieties of small codimensions|Proc.\ of Symp.\ in Pure
Math.|29|1975|129--164|

\livre Hi66|Hirzebruch F|Topological methods in Algebraic Geometry|
Grundl.\ Math.\ Wiss.{\bf 131}, Spriger, Heidelberg|1966|

\article Jo78|Jouanolou J.-P|Hypersurfaces solutions d'une \'equation de Pfaff
analytique| Math.\ Ann.|332|1978|239--248| 

\livre Ko70|Kobayashi S|Hyperbolic manifolds and holomorphic
mappings|Marcel Dek\-ker, New York|1970|

\article La86|Lang S|Hyperbolic and diophantine
analysis|Bull.\ Amer.\ Math.\ Soc.|14|1986|159--205|

\livre La87|Lang S|Introduction to complex hyperbolic
spaces|Springer-Verlag, New York|1987|

\article LY90|Lu S.S.-Y., Yau S.T|Holomorphic curves in surfaces of
general type|Proc.\ Nat.\ Acad.\ Sci.\ USA|87|1990|80--82|

\article M-De78|Martin-Deschamps M|Courbes de genre g\'eom\'etrique born\'e
sur une surface de type g\'en\'eral (d'apr\`es F.~Bogomolov)|S\'eminaire
Bourbaki|519|1978|1--15|

\divers McQ97|McQuillan M| Diophantine approximations and foliations|
Preprint (1997), to appear in Publ.\ Math.\ IHES| 

\article Mi82|Miyaoka Y|Algebraic surfaces with positive
indices|Katata Symp.\ Proc.\ 1982, Progress in Math.|39|1983|281--301|

\divers Pa99|Paun, M|Personal communication|January 1999|

\livre Sa78|Sakai F|Symmetric powers of the cotangent bundle and
classification of algebraic varieties|Proc.\ Copenhagen Meeting in
Alg. Geom.|1978|

\article ST88|Schneider M., Tancredi A|Almost-positive vector bundles
on projective surfaces|Math.\ Ann.|280|1988|537--547|

\article Se68|Seidenberg A| Reduction of singularities of the differential
equation $AdY=BdX$|Amer.\ J.\ of Math.|90|1968|248--269|

\article Siu87|Siu Y.-T|Defect relations for holomorphic maps between
spaces of different dimensions|Duke Math.\ J.|55|1987|213--251|

\livre SY95|Siu Y.-T., Yeung S.K|Hyperbolicity of the complement of a
generic smooth curve of high degree in the complex projective plane| to
appear in Inventiones Math.|1996|

\divers SY97|Siu Y.-T., Yeung S.K|Defects for ample divisors of abelian
varieties, Schwarz lemma and hyperbolic hypersurfaces of low
degrees|Preprint, Fall 1996|

\article Voi96|Voisin, C|On a conjecture of Clemens on rational curves 
on hypersurfaces|J.\ Differ.\ Geom|44|1996|200--213|

\article Xu94|Xu G|Subvarieties of general hypersurfaces in projective
space|J.\ Differential Geometry|39|1994|139--172|

\article Yau78|Yau S.T|On the Ricci curvature of a compact K\"ahler manifold
and the complex Monge Amp\`ere equation|Comm.\ Pure and Appl.\ Math.|31|1978|
339--441|

\article Zai87|Zaidenberg, M|The complement of a generic
hypersurface of degree $2n$ in $\bC\bP^n$ is not hyperbolic|
Siberian Math.\ J.|28|1987|425--432|

\article Zai89|Zaidenberg, M|Stability of hyperbolic embeddedness
and construction of examples|Math.\ USSR Sbornik|63|1989|351--361|
 
\divers Zai93|Zaidenberg, M|Hyperbolicity in projective spaces|International
Symposium on Holomorphic mappings, Diophantine Geometry
and Related topics, R.I.M.S.\ Lecture Notes ser.\ {\bf 819},
R.I.M.S.\ Kyoto University (1993), 136--156|

}

\parindent=0cm
\vskip10pt
(printed on \today, \timeofday)
\vskip20pt

Universit\'e Joseph Fourier Grenoble I\hfil\break
Institut Fourier (Math\'ematiques)\hfil\break
UMR 5582 du C.N.R.S., BP 74\hfil\break 38402
Saint-Martin d'H\`eres, France\hfil\break
{\it e-mail:}\/ demailly@fourier.ujf-grenoble.fr, 
elgoul@picard.ups-tlse.fr

\immediate\closeout1
\end